\ifx\shlhetal\undefinedcontrolsequence\let\shlhetal\relax\fi

%
\input amstex
\expandafter\ifx\csname mathdefs.tex\endcsname\relax
  \expandafter\gdef\csname mathdefs.tex\endcsname{}
\else \message{Hey!  Apparently you were trying to
  \string twice.   This does not make sense.} 
\errmessage{Please edit your file (probably \jobname.tex) and remove
any duplicate ``\string\input'' lines} \fi




\catcode`\X=12\catcode`\@=11

\def\n@wcount{\alloc@0\count\countdef\insc@unt}
\def\n@wwrite{\alloc@7\write\chardef\sixt@@n}
\def\n@wread{\alloc@6\read\chardef\sixt@@n}
\def\r@s@t{\relax}\def\v@idline{\par}\def\@mputate#1/{#1}
\def\l@c@l#1X{\firstpart.#1}\def\gl@b@l#1X{#1}\def\t@d@l#1X{{}}

\def\crossrefs#1{\ifx\all#1\let\tr@ce=\all\else\def\tr@ce{#1,}\fi
   \n@wwrite\cit@tionsout\openout\cit@tionsout=\jobname.cit 
   \write\cit@tionsout{\tr@ce}\expandafter\setfl@gs\tr@ce,}
\def\setfl@gs#1,{\def\@{#1}\ifx\@\empty\let\next=\relax
   \else\let\next=\setfl@gs\expandafter\xdef
   \csname#1tr@cetrue\endcsname{}\fi\next}
\def\m@ketag#1#2{\expandafter\n@wcount\csname#2tagno\endcsname
     \csname#2tagno\endcsname=0\let\tail=\all\xdef\all{\tail#2,}
   \ifx#1\l@c@l\let\tail=\r@s@t\xdef\r@s@t{\csname#2tagno\endcsname=0\tail}\fi
   \expandafter\gdef\csname#2cite\endcsname##1{\expandafter
     \ifx\csname#2tag##1\endcsname\relax?\else\csname#2tag##1\endcsname\fi
     \expandafter\ifx\csname#2tr@cetrue\endcsname\relax\else
     \write\cit@tionsout{#2tag ##1 cited on page \folio.}\fi}
   \expandafter\gdef\csname#2page\endcsname##1{\expandafter
     \ifx\csname#2page##1\endcsname\relax?\else\csname#2page##1\endcsname\fi
     \expandafter\ifx\csname#2tr@cetrue\endcsname\relax\else
     \write\cit@tionsout{#2tag ##1 cited on page \folio.}\fi}
   \expandafter\gdef\csname#2tag\endcsname##1{\expandafter
      \ifx\csname#2check##1\endcsname\relax
      \expandafter\xdef\csname#2check##1\endcsname{}%
      \else\immediate\write16{Warning: #2tag ##1 used more than once.}\fi
      \multit@g{#1}{#2}##1/X%
      \write\t@gsout{#2tag ##1 assigned number \csname#2tag##1\endcsname\space
      on page \number\count0.}%
   \csname#2tag##1\endcsname}}
\def\multit@g#1#2#3/#4X{\def\t@mp{#4}\ifx\t@mp\empty%
      \global\advance\csname#2tagno\endcsname by 1 
      \expandafter\xdef\csname#2tag#3\endcsname
      {#1\number\csname#2tagno\endcsnameX}%
   \else\expandafter\ifx\csname#2last#3\endcsname\relax
      \expandafter\n@wcount\csname#2last#3\endcsname
      \global\advance\csname#2tagno\endcsname by 1 
      \expandafter\xdef\csname#2tag#3\endcsname
      {#1\number\csname#2tagno\endcsnameX}
      \write\t@gsout{#2tag #3 assigned number \csname#2tag#3\endcsname\space
      on page \number\count0.}\fi
   \global\advance\csname#2last#3\endcsname by 1
   \def\t@mp{\expandafter\xdef\csname#2tag#3/}%
   \expandafter\t@mp\@mputate#4\endcsname
   {\csname#2tag#3\endcsname\lastpart{\csname#2last#3\endcsname}}\fi}
\def\t@gs#1{\def\all{}\m@ketag#1e\m@ketag#1s\m@ketag\t@d@l p
   \m@ketag\gl@b@l r \n@wread\t@gsin
   \openin\t@gsin=\jobname.tgs \re@der \closein\t@gsin
   \n@wwrite\t@gsout\openout\t@gsout=\jobname.tgs }
\outer\def\localtags{\t@gs\l@c@l}
\outer\def\globaltags{\t@gs\gl@b@l}
\outer\def\newlocaltag#1{\m@ketag\l@c@l{#1}}
\outer\def\newglobaltag#1{\m@ketag\gl@b@l{#1}}

\newif\ifpr@ 
\def\m@kecs #1tag #2 assigned number #3 on page #4.%
   {\expandafter\gdef\csname#1tag#2\endcsname{#3}
   \expandafter\gdef\csname#1page#2\endcsname{#4}
   \ifpr@\expandafter\xdef\csname#1check#2\endcsname{}\fi}
\def\re@der{\ifeof\t@gsin\let\next=\relax\else
   \read\t@gsin to\t@gline\ifx\t@gline\v@idline\else
   \expandafter\m@kecs \t@gline\fi\let \next=\re@der\fi\next}
\def\pretags#1{\pr@true\pret@gs#1,,}
\def\pret@gs#1,{\def\@{#1}\ifx\@\empty\let\n@xtfile=\relax
   \else\let\n@xtfile=\pret@gs \openin\t@gsin=#1.tgs \message{#1} \re@der 
   \closein\t@gsin\fi \n@xtfile}

\newcount\sectno\sectno=0\newcount\subsectno\subsectno=0
\newif\ifultr@local \def\ultralocal{\ultr@localtrue}
\def\firstpart{\number\sectno}
\def\lastpart#1{\ifcase#1 \or a\or b\or c\or d\or e\or f\or g\or h\or 
   i\or k\or l\or m\or n\or o\or p\or q\or r\or s\or t\or u\or v\or w\or 
   x\or y\or z \fi}

\def\resetall{\global\advance\sectno by 1\subsectno=0
   \gdef\firstpart{\number\sectno}\r@s@t}
\def\resetsub{\global\advance\subsectno by 1
   \gdef\firstpart{\number\sectno.\number\subsectno}\r@s@t}
\def\newsection#1\par{\resetall\vskip0pt plus.3\vsize\penalty-250
   \vskip0pt plus-.3\vsize\bigskip\bigskip
   \message{#1}\leftline{\bf#1}\nobreak\bigskip}
\def\subsection#1\par{\ifultr@local\resetsub\fi
   \vskip0pt plus.2\vsize\penalty-250\vskip0pt plus-.2\vsize
   \bigskip\smallskip\message{#1}\leftline{\bf#1}\nobreak\medskip}

\def\t@gsoff#1,{\def\@{#1}\ifx\@\empty\let\next=\relax\else\let\next=\t@gsoff
   \def\@@{p}\ifx\@\@@\else
   \expandafter\gdef\csname#1cite\endcsname##1{\zeigen{##1}}
   \expandafter\gdef\csname#1page\endcsname##1{?}
   \expandafter\gdef\csname#1tag\endcsname##1{\zeigen{##1}}\fi\fi\next}
\def\verbatimtags{\ifx\all\relax\else\expandafter\t@gsoff\all,\fi}
\def\zeigen#1{\hbox{$\langle$}#1\hbox{$\rangle$}}

\def\(#1){\edef\dot@g{\ifmmode\ifinner(\hbox{\noexpand\etag{#1}})
   \else\noexpand\eqno(\hbox{\noexpand\etag{#1}})\fi
   \else(\noexpand\ecite{#1})\fi}\dot@g}

\newif\ifbr@ck
\def\eat#1{}
\def\[#1]{\br@cktrue[\br@cket#1'X]}
\def\br@cket#1'#2X{\def\temp{#2}\ifx\temp\empty\let\next\eat
   \else\let\next\br@cket\fi
   \ifbr@ck\br@ckfalse\br@ck@t#1,X\else\br@cktrue#1\fi\next#2X}
\def\br@ck@t#1,#2X{\def\temp{#2}\ifx\temp\empty\let\neext\eat
   \else\let\neext\br@ck@t\def\temp{,}\fi
   \def\teemp{#1}\ifx\teemp\empty\else\rcite{#1}\fi\temp\neext#2X}
\def\resetbr@cket{\gdef\[##1]{[\rtag{##1}]}}
\def\references{\resetbr@cket\newsection References\par}

\newtoks\symb@ls\newtoks\s@mb@ls\newtoks\p@gelist\n@wcount\ftn@mber
    \ftn@mber=1\newif\ifftn@mbers\ftn@mbersfalse\newif\ifbyp@ge\byp@gefalse
\def\defm@rk{\ifftn@mbers\n@mberm@rk\else\symb@lm@rk\fi}
\def\n@mberm@rk{\xdef\m@rk{{\the\ftn@mber}}%
    \global\advance\ftn@mber by 1 }
\def\rot@te#1{\let\temp=#1\global#1=\expandafter\r@t@te\the\temp,X}
\def\r@t@te#1,#2X{{#2#1}\xdef\m@rk{{#1}}}
\def\b@@st#1{{$^{#1}$}}\def\str@p#1{#1}
\def\symb@lm@rk{\ifbyp@ge\rot@te\p@gelist\ifnum\expandafter\str@p\m@rk=1 
    \s@mb@ls=\symb@ls\fi\write\f@nsout{\number\count0}\fi \rot@te\s@mb@ls}
\def\byp@ge{\byp@getrue\n@wwrite\f@nsin\openin\f@nsin=\jobname.fns 
    \n@wcount\currentp@ge\currentp@ge=0\p@gelist={0}
    \re@dfns\closein\f@nsin\rot@te\p@gelist
    \n@wread\f@nsout\openout\f@nsout=\jobname.fns }
\def\m@kelist#1X#2{{#1,#2}}
\def\re@dfns{\ifeof\f@nsin\let\next=\relax\else\read\f@nsin to \f@nline
    \ifx\f@nline\v@idline\else\let\t@mplist=\p@gelist
    \ifnum\currentp@ge=\f@nline
    \global\p@gelist=\expandafter\m@kelist\the\t@mplistX0
    \else\currentp@ge=\f@nline
    \global\p@gelist=\expandafter\m@kelist\the\t@mplistX1\fi\fi
    \let\next=\re@dfns\fi\next}
\def\symbols#1{\symb@ls={#1}\s@mb@ls=\symb@ls} 
\def\bigsymbol{\textstyle}
\symbols{\bigsymbol\ast,\dagger,\ddagger,\sharp,\flat,\natural,\star}
\def\ftnumbers{\ftn@mberstrue} \def\ftsymbols{\ftn@mbersfalse}
\def\paginal{\byp@ge} \def\resetftnumbers{\ftn@mber=1}
\def\ftnote#1{\defm@rk\expandafter\expandafter\expandafter\footnote
    \expandafter\b@@st\m@rk{#1}}

\long\def\jump#1\endjump{}
\def\ssum{\mathop{\lower .1em\hbox{$\textstyle\Sigma$}}\nolimits}

\def\qed{\nobreak\kern 1em \vrule height .5em width .5em depth 0em}
\def\newneq{\hbox{\rlap{\hbox to 1\wd9{\hss$=$\hss}}\raise .1em 
   \hbox to 1\wd9{\hss$\scriptscriptstyle/$\hss}}}
\def\subsetne{\setbox9 = \hbox{$\subset$}\mathrel{\hbox{\rlap
   {\lower .4em \newneq}\raise .13em \hbox{$\subset$}}}}
\def\supsetne{\setbox9 = \hbox{$\subset$}\mathrel{\hbox{\rlap
   {\lower .4em \newneq}\raise .13em \hbox{$\supset$}}}}

\def\vbar{\mathchoice{\vrule height6.3ptdepth-.5ptwidth.8pt\kern-.8pt}
   {\vrule height6.3ptdepth-.5ptwidth.8pt\kern-.8pt}
   {\vrule height4.1ptdepth-.35ptwidth.6pt\kern-.6pt}
   {\vrule height3.1ptdepth-.25ptwidth.5pt\kern-.5pt}}
\def\f@dge{\mathchoice{}{}{\mkern.5mu}{\mkern.8mu}}
\def\b@c#1#2{{\rm \mkern#2mu\vbar\mkern-#2mu#1}}
\def\b@b#1{{\rm I\mkern-3.5mu #1}}
\def\b@a#1#2{{\rm #1\mkern-#2mu\f@dge #1}}
\def\bb#1{{\count4=`#1 \advance\count4by-64 \ifcase\count4\or\b@a A{11.5}\or
   \b@b B\or\b@c C{5}\or\b@b D\or\b@b E\or\b@b F \or\b@c G{5}\or\b@b H\or
   \b@b I\or\b@c J{3}\or\b@b K\or\b@b L \or\b@b M\or\b@b N\or\b@c O{5} \or
   \b@b P\or\b@c Q{5}\or\b@b R\or\b@a S{8}\or\b@a T{10.5}\or\b@c U{5}\or
   \b@a V{12}\or\b@a W{16.5}\or\b@a X{11}\or\b@a Y{11.7}\or\b@a Z{7.5}\fi}}

\catcode`\X=11 \catcode`\@=12

\sectno=0   
\localtags
\NoBlackBoxes
\def\cite #1{\rm[#1]}
\def\lesscircle{{\mathrel{\mathord{<}
  \!\!
  \raise 0.8 pt\hbox{$\scriptstyle\circ$}}}}
\documentstyle {amsppt}
\topmatter
\title {Historic iteration with $\aleph_\varepsilon$-support \\
Sh538} \endtitle
\author {Saharon Shelah \thanks {\null\newline
I thank Otmar Spinas, Juris Steprans and Boban Velickovic for the
helpful discussion. \null\newline
I thank Alice Leonhardt for the beautiful typing. \null\newline
Latest Revision 96/Aug/2 \null\newline
Originally typed 94/Jan/23 elsewhere} \endthanks} \endauthor
\bigskip
\affil {Institute of Mathematics \\
The Hebrew University \\
Jerusalem, Israel
\medskip
Rutgers University \\
Department of Mathematics \\
New Brunswick, NJ USA} \endaffil
\bigskip

\abstract  One aim of this work is to get a universe in which weak
versions of Martin axioms holds for some forcing notions of cardinality
$\aleph_0,\aleph_1$ and $\aleph_2$ while on $\aleph_2$ club, the ``small"
brother of diamond, holds.  
As a consequence we get the consistency of ``there is
no Gross space".  Another aim is to present a case of ``Historic iteration".
\endabstract
\endtopmatter
\document  
\expandafter\ifx\csname bib4plain.tex\endcsname\relax
  \expandafter\gdef\csname bib4plain.tex\endcsname{}
\else \message{Hey!  Apparently you were trying to \string twice.   This does not make sense.}
\errmessage{Please edit your file (probably \jobname.tex) and remove
any duplicate ``\string\input'' lines} \fi

\def\renewcommand{\newcommand}	       
\edef\cite{\the\catcode`@}%
\catcode`@ = 11
\let\@oldatcatcode = \cite
\chardef\@letter = 11
\chardef\@other = 12
%
%
%
%
\def\@innerdef#1#2{\edef#1{\expandafter\noexpand\csname #2\endcsname}}%
%
%
\@innerdef\@innernewcount{newcount}%
\@innerdef\@innernewdimen{newdimen}%
\@innerdef\@innernewif{newif}%
\@innerdef\@innernewwrite{newwrite}%
%
%
%
\def\@gobble#1{}%
%
%
%
\ifx\inputlineno\@undefined
   \let\@linenumber = \empty 
\else
   \def\@linenumber{\the\inputlineno:\space}%
\fi
%
%
%
\def\@futurenonspacelet#1{\def\cs{#1}%
   \afterassignment\@stepone\let\@nexttoken=
}%
\begingroup 
\def\\{\global\let\@stoken= }%
\\ 
\endgroup
\def\@stepone{\expandafter\futurelet\cs\@steptwo}%
\def\@steptwo{\expandafter\ifx\cs\@stoken\let\@@next=\@stepthree
   \else\let\@@next=\@nexttoken\fi \@@next}%
\def\@stepthree{\afterassignment\@stepone\let\@@next= }%
%
%
%
\def\@getoptionalarg#1{%
   \let\@optionaltemp = #1%
   \let\@optionalnext = \relax
   \@futurenonspacelet\@optionalnext\@bracketcheck
}%
%
%
\def\@bracketcheck{%
   \ifx [\@optionalnext
      \expandafter\@@getoptionalarg
   \else
      \let\@optionalarg = \empty
      \expandafter\@optionaltemp
   \fi
}%
\def\@@getoptionalarg[#1]{%
   \def\@optionalarg{#1}%
   \@optionaltemp
}%
%
%
%
\def\@nnil{\@nil}%
\def\@fornoop#1\@@#2#3{}%
\def\@for#1:=#2\do#3{%
   \edef\@fortmp{#2}%
   \ifx\@fortmp\empty \else
      \expandafter\@forloop#2,\@nil,\@nil\@@#1{#3}%
   \fi
}%
\def\@forloop#1,#2,#3\@@#4#5{\def#4{#1}\ifx #4\@nnil \else
       #5\def#4{#2}\ifx #4\@nnil \else#5\@iforloop #3\@@#4{#5}\fi\fi
}%
\def\@iforloop#1,#2\@@#3#4{\def#3{#1}\ifx #3\@nnil
       \let\@nextwhile=\@fornoop \else
      #4\relax\let\@nextwhile=\@iforloop\fi\@nextwhile#2\@@#3{#4}%
}%
%
%
%
\@innernewif\if@fileexists
\def\@testfileexistence{\@getoptionalarg\@finishtestfileexistence}%
\def\@finishtestfileexistence#1{%
   \begingroup
      \def\extension{#1}%
      \immediate\openin0 =
         \ifx\@optionalarg\empty\jobname\else\@optionalarg\fi
         \ifx\extension\empty \else .#1\fi
         \space
      \ifeof 0
         \global\@fileexistsfalse
      \else
         \global\@fileexiststrue
      \fi
      \immediate\closein0
   \endgroup
}%
%
%
%
%
\def\bibliographystyle#1{%
   \@readauxfile
   \@writeaux{\string\bibstyle{#1}}%
}%
\let\bibstyle = \@gobble
%
%
\let\bblfilebasename = \jobname
\def\bibliography#1{%
   \@readauxfile
   \@writeaux{\string\bibdata{#1}}%
   \@testfileexistence[\bblfilebasename]{bbl}%
   \if@fileexists
      \nobreak
      \@readbblfile
   \fi
}%
\let\bibdata = \@gobble
%
%
\def\nocite#1{%
   \@readauxfile
   \@writeaux{\string\citation{#1}}%
}%
\@innernewif\if@notfirstcitation
%
%
\def\cite{\@getoptionalarg\@cite}%
%
%
\def\@cite#1{%
   \let\@citenotetext = \@optionalarg
   \printcitestart
   \nocite{#1}%
   \@notfirstcitationfalse
   \@for \@citation :=#1\do
   {%
      \expandafter\@onecitation\@citation\@@
   }%
   \ifx\empty\@citenotetext\else
      \printcitenote{\@citenotetext}%
   \fi
   \printcitefinish
}%
\def\@onecitation#1\@@{%
   \if@notfirstcitation
      \printbetweencitations
   \fi
   \expandafter \ifx \csname\@citelabel{#1}\endcsname \relax
      \if@citewarning
         \message{\@linenumber Undefined citation `#1'.}%
      \fi
      \expandafter\gdef\csname\@citelabel{#1}\endcsname{%
\strut
\vadjust{\vskip-\dp\strutbox
\vbox to 0pt{\vss\parindent0cm \leftskip=\hsize 
\advance\leftskip3mm
\advance\hsize 4cm\strut\openup-4pt 
\rightskip 0cm plus 1cm minus 0.5cm ?  #1 ?\strut}}
         {\tt
            \escapechar = -1
            \nobreak\hskip0pt
            \expandafter\string\csname#1\endcsname
            \nobreak\hskip0pt
         }%
      }%
   \fi
   \csname\@citelabel{#1}\endcsname
   \@notfirstcitationtrue
}%
%
%
\def\@citelabel#1{b@#1}%
%
%
\def\@citedef#1#2{\expandafter\gdef\csname\@citelabel{#1}\endcsname{#2}}%
%
%
%
\def\@readbblfile{%
   \ifx\@itemnum\@undefined
      \@innernewcount\@itemnum
   \fi
   \begingroup
      \def\begin##1##2{%
         \setbox0 = \hbox{\biblabelcontents{##2}}%
         \biblabelwidth = \wd0
      }%
      \def\end##1{}
      %
      %
      \@itemnum = 0
      \def\bibitem{\@getoptionalarg\@bibitem}%
      \def\@bibitem{%
         \ifx\@optionalarg\empty
            \expandafter\@numberedbibitem
         \else
            \expandafter\@alphabibitem
         \fi
      }%
      \def\@alphabibitem##1{%
         \expandafter \xdef\csname\@citelabel{##1}\endcsname {\@optionalarg}%
         \ifx\biblabelprecontents\@undefined
            \let\biblabelprecontents = \relax
         \fi
         \ifx\biblabelpostcontents\@undefined
            \let\biblabelpostcontents = \hss
         \fi
         \@finishbibitem{##1}%
      }%
      \def\@numberedbibitem##1{%
         \advance\@itemnum by 1
         \expandafter \xdef\csname\@citelabel{##1}\endcsname{\number\@itemnum}%
         \ifx\biblabelprecontents\@undefined
            \let\biblabelprecontents = \hss
         \fi
         \ifx\biblabelpostcontents\@undefined
            \let\biblabelpostcontents = \relax
         \fi
         \@finishbibitem{##1}%
      }%
      \def\@finishbibitem##1{%
         \biblabelprint{\csname\@citelabel{##1}\endcsname}%
         \@writeaux{\string\@citedef{##1}{\csname\@citelabel{##1}\endcsname}}%
         \ignorespaces
      }%
      %
      %
      \let\em = \bblem
      \let\newblock = \bblnewblock
      \let\sc = \bblsc
      \frenchspacing
      \clubpenalty = 4000 \widowpenalty = 4000
      \tolerance = 10000 \hfuzz = .5pt
      \everypar = {\hangindent = \biblabelwidth
                      \advance\hangindent by \biblabelextraspace}%
      \bblrm
      \parskip = 1.5ex plus .5ex minus .5ex
      \biblabelextraspace = .5em
      \bblhook
      \input \bblfilebasename.bbl
   \endgroup
}%
%
%
\@innernewdimen\biblabelwidth
\@innernewdimen\biblabelextraspace
%
%
%
\def\biblabelprint#1{%
   \noindent
   \hbox to \biblabelwidth{%
      \biblabelprecontents
      \biblabelcontents{#1}%
      \biblabelpostcontents
   }%
   \kern\biblabelextraspace
}%
%
%
%
\def\biblabelcontents#1{{\bblrm [#1]}}%
%
%
\def\bblrm{\rm}%
%
%
\def\bblem{\it}%
%
%
\def\bblsc{\ifx\@scfont\@undefined
              \font\@scfont = cmcsc10
           \fi
           \@scfont
}%
%
%
\def\bblnewblock{\hskip .11em plus .33em minus .07em }%
%
%
\let\bblhook = \empty
%
%
%
\def\printcitestart{[}
\def\printcitefinish{]}
\def\printbetweencitations{, }
\def\printcitenote#1{, #1}
%
%
%
\let\citation = \@gobble
%
%
%
\@innernewcount\@numparams
%
%
\def\newcommand#1{%
   \def\@commandname{#1}%
   \@getoptionalarg\@continuenewcommand
}%
%
%
\def\@continuenewcommand{%
   \@numparams = \ifx\@optionalarg\empty 0\else\@optionalarg \fi \relax
   \@newcommand
}%
%
%
\def\@newcommand#1{%
   \def\@startdef{\expandafter\edef\@commandname}%
   \ifnum\@numparams=0
      \let\@paramdef = \empty
   \else
      \ifnum\@numparams>9
         \errmessage{\the\@numparams\space is too many parameters}%
      \else
         \ifnum\@numparams<0
            \errmessage{\the\@numparams\space is too few parameters}%
         \else
            \edef\@paramdef{%
               \ifcase\@numparams
                  \empty  No arguments.
               \or ####1%
               \or ####1####2%
               \or ####1####2####3%
               \or ####1####2####3####4%
               \or ####1####2####3####4####5%
               \or ####1####2####3####4####5####6%
               \or ####1####2####3####4####5####6####7%
               \or ####1####2####3####4####5####6####7####8%
               \or ####1####2####3####4####5####6####7####8####9%
               \fi
            }%
         \fi
      \fi
   \fi
   \expandafter\@startdef\@paramdef{#1}%
}%
%
%
%
%
\def\@readauxfile{%
   \if@auxfiledone \else 
      \global\@auxfiledonetrue
      \@testfileexistence{aux}%
      \if@fileexists
         \begingroup
            \endlinechar = -1
            \catcode`@ = 11
            \input \jobname.aux
         \endgroup
      \else
         \message{\@undefinedmessage}%
         \global\@citewarningfalse
      \fi
      \immediate\openout\@auxfile = \jobname.aux
   \fi
}%
%
%
\newif\if@auxfiledone
\ifx\noauxfile\@undefined \else \@auxfiledonetrue\fi
%
%
%
%
\@innernewwrite\@auxfile
\def\@writeaux#1{\ifx\noauxfile\@undefined \write\@auxfile{#1}\fi}%
%
%
%
\ifx\@undefinedmessage\@undefined
   \def\@undefinedmessage{No .aux file; I won't give you warnings about
                          undefined citations.}%
\fi
%
%
\@innernewif\if@citewarning
\ifx\noauxfile\@undefined \@citewarningtrue\fi
%
%
%
\catcode`@ = \@oldatcatcode
 
\newpage

\head {\S1} \endhead
\resetall
\bigskip

When we make the continuum $> \aleph_2$, it in general takes effort not to add
$\aleph_2$ Cohen reals.  We do it here, making a weak version of MA true 
without adding $\aleph_2$ Cohen reals (i.e. a function $f:\omega_2 
\rightarrow \{0,1\}$ which is quite generic for the forcing notion
$\{g:g$ a finite function from $\omega_2$ to $\{0,1\}\}$).  Moreover 
$\clubsuit_{S^2_0}$ is
preserved.  This was motivated by the following application to Gross spaces:
the consistency of ``there is no Gross space" with ZFC, prove in \scite{1.23} 
below; on Gross spaces see e.g. Shelah Spinas \cite{ShSi:468}.

We prove that for a restricted enough family of $\sigma$-centered forcing,
(essentially giving almost intersection to ``definable" filters) we have:
\medskip
\roster
\item "{$(a)$}"  for any $\aleph_1$ condtions, $\aleph_1$ of them belongs to a
directed subset
\item "{$(b)$}"  for any 
$\aleph_2$ conditions, some $\aleph_2$ of them belongs to a directed subset
\item "{$(c)$}"  $\clubsuit_{\aleph_2}$
\item "{$(d)$}"  $2^{\aleph_0} = \aleph_3$.
\endroster
\medskip

\noindent
This is enough; we can by the proof strengthen (a) (even to
MA$_{\aleph_1}(\sigma$-centered)), it is less clear about (b).  This was our
first approach, but it seems considerably clearer to prove just (b).
\bigskip

Definition \scite{1.1} defines the family of c.c.c. forcing notions 
for which we will have our approximation to (weak) MA.  A main Definition is 
\scite{1.11}, where we
define the family of iterations we shall use.  For preserving $\clubsuit
_{S^2_0}$ we want, among any given $\aleph_2$ conditions $\langle p_i:i <
\omega_2 \rangle$ to find (quite many) countable subsets $w \subseteq
\omega_2$ such that $\{p_i:i \in w \}$ has an upper bound.  So we have 
to ``marry" the c.c.c. with somewhat countable support hence the name
$\aleph_\varepsilon$-support.  The support is fashioned
after ``historic forcing" (see Shelah Stanley \cite{ShSt:258}).  
We may continue this generalization to a larger family of forcing, 
replacing $\aleph_0$ by a larger cardinal (see \scite{1.3}(5)).

Note: $\clubsuit_{S^2_0} \Rightarrow$ for some $A \subseteq \aleph_2$ there
is no $L[A]$-generic $\omega_2$-sequence of Cohen reals $\Rightarrow {\frak p}
\le \aleph_2$ (where ${\frak p} \le \aleph_2$ means for some $A_i \in D$
(for $i < \omega_2)$, $D$ is a filter on $\omega$ containing all co-finite 
sets, there is no $A \in [\omega]^{\aleph_0}$ such that $(\forall i < 
\omega_2)[A \subseteq^* A_i]$).

In \scite{1.1}-\scite{1.2} we deal with simple forcing, a very 
restricted class but sufficient for our need.  In \scite{1.3} we comment on
generalizations, in \scite{1.4} we fix some cardinal parameters.
In \scite{1.5} - \scite{1.10} we deal with ``creatures", they 
represent possible ``support" of a condition in the 
iteration to be defined later.  In \scite{1.11} we define our iteration and
prove basic facts by simultaneous induction of length.  In \scite{1.12} -
\scite{1.17} we further investigate the iteration.  In \scite{1.18} -
\scite{1.22} we get the actual consistency results.
\bigskip

\noindent
\underbar{Notation}:  Let $u,v$ denote creatures (see Definition \scite{1.5} 
below), $\bold p, \bold q, \bold r$ denote conditions, (in the iteration), 
$\tau,\xi,\zeta,
\varepsilon$ denote countable ordinals (for depth of creatures or conditions)
and $\alpha,\beta,\gamma,\delta$ denote ordinals.  We use $\varphi,\psi,
\Phi,\Psi$ only as in Definition \scite{1.5}.
\bigskip

\noindent
The following looks extremely special but 
includes the forcing notion needed for the
application to Gross spaces.  See Example \scite{1.3}(3); we may consider some
generalizations.
\definition{\stag{1.1} Definition}  1) A forcing notion $Q$ is a 
simple forcing notion if:
\medskip
\roster
\item "{$(a)$}"  for some $S \subseteq \dsize \bigcup_{n < \omega}
{}^n \omega:p \in Q$ \underbar{iff} $p$ has the form $(s,\bold t)$ where:
\newline
$s \in S$ and $\bold t$ is a finite subset of ${}^\omega \omega$
\item "{$(b)$}"  $(s^1,\bold t^1) \le (s^2,\bold t^2)$ implies $s^1 \subseteq
s^2 \and \bold t^1 \subseteq \bold t^2$
\item "{$(c)$}"  $(s,\bold t') \le (s,\bold t'')$ whenever $\bold t'
\subseteq \bold t''$ and $s \in S$, also \newline
$(\emptyset,\emptyset) \in Q$ and
$(\emptyset,\emptyset) \le (s,\bold t)$
\item "{$(d)$}"   for some two-place functions $f,g$ 
[with domain and range $\subseteq {\Cal H}(\aleph_0)$ - \newline
called the witnesses], \underbar{if} $(s,\bold t) \in Q,m <
\omega,|\{t \restriction k:t \in \bold t\}| \le m$, \newline
$\text{Dom}(s) = n
\text{ and } k=f(n,m)$ \underbar{then} $(s,\bold t) \le (s^+,\bold t) \in Q$
\newline
and $s^+(n) \ne 0$ where 
$s^+ = s \bigcup \{ \langle n,g(s,\{t_\ell \restriction k:\ell < m\})
\rangle\}$ (so part of the conclusion is
``$g(s,\{t_\ell \restriction k:\ell < m\})$ is well defined"; we allow $f$ and
$g$ to be defined in additional cases for technical reasons)
\item "{$(e)$}"  if $(s,\bold t) \in Q,\text{ Dom}(s) = n$ then $(s,\bold t)
\le (s \bigcup \{ \langle n,0 \rangle\},\bold t) \in Q$
\item "{$(f)$}"  the truth value of $(s^1_2,\bold t^1) \le (s^2_2,\bold t^2)$
depend just on \newline
$(s_1,\{\eta \restriction \ell g(s^2):\eta \in \bold t^1\},
s_2,\{\eta \restriction \ell g(s^2):\eta \in \bold t^2\})$ \newline
(not really needed but natural). 
\endroster
\medskip

\noindent
2) If $\bold t^1,\bold t^2 \subseteq {}^\omega \omega$ let
$(s^1,\bold t^1) \le_Q (s^2,\bold t^2)$ means: for every finite
$\bold t_1$ of $\bold t^1$ for some finite subset $\bold t_2$ of $\bold t^2$
we have $(s^1,\bold t_1) \le_Q (s^2,\bold t_2)$. \newline
3)  If $Q$ is a simple forcing notion we say $X$ describes $Q$ if
$X = S_Q = (S,Y)$ where $Y =: \{(s,n,\ell,k,\bar t):s \in S, t = \langle
t_\ell:\ell < \ell^* \rangle, t_\ell \in {}^k \omega$ and: $\bold t \subseteq
{}^\omega \omega$ finite, $\{t \restriction k:t \in \bold t\} \subseteq
\{t_\ell:\ell < \ell^*\} \Rightarrow (s,\bold t) \le (s \char 94 \langle \ell
\rangle,\bold t)\}$ (so from $X$, we can compute $f,g$ (more exactly, the
Borel set of pairs $(f,g)$ which are as required).)
\bigskip

\noindent
\underbar{\stag{1.2} Examples}:  The following are simple forcing notions:
\newline
1)  Cohen forcing, here $S \equiv \dsize \bigcup_{n < \omega} {}^n 2$ and
$(s^1,\bold t^1) \le (s^2,\bold t^2)$ iff $s^1 \subseteq
s^2,\bold t^1 \subseteq \bold t^2$.  Here \newline
$f(n,m) = 1, g(s,\{t_\ell
\restriction k:\ell < m\}) = 1$ are O.K. \newline
2)  Dominating real forcing $=$ Hechler forcing: $(s^1,\bold t^1) \le
(s^2,\bold t^2)$ \underbar{iff} $s^1 \subseteq s^2,\bold t^1 \subseteq
\bold t^2$ and $[n \in \text{ Dom}(s^2) \backslash \text{ Dom}(s^1) \and
t \in \bold t^1 \and s^2(n) \ne 0 \Rightarrow s^2(n) \ge t(n)]$.  Here 
$f(n,m) = n+1, g(s,\{t_\ell \restriction k:\ell < m\}) = \text{ Max}
[\{(t_\ell \restriction k)(n) + 1:\ell < m\} \bigcup \{1\}]$ are O.K. \newline
3)   Vector spaces forcing $=$ Spinas forcing: let $K$ be a countable field,
${\Cal U}$ a vector space with dimension $\aleph_0$ and free basis 
$X = \{x_n:n < \omega \}$; without loss of generality $K \subseteq K^*,K^*$ 
is a countable algebraically closed field with transcendence dimension
$\aleph_0$ and $0_{K^*}$ (the zero of $K^*$) is the ordinal $0,{\Cal U}^*$ 
is the $K^*$-vector space extending ${\Cal U}$ and they have a common basic 
$X$ and $0_{\Cal U}^* = 0$ and the set of elements of $K^*$ is $\omega$ and 
also the set of elements of ${\Cal U}^*$ is exactly $\omega$.  We can assume
that $K^*,{\Cal U}^* \in L$. For $f:\omega \rightarrow K^*$ let 
$f^{\text{HOM}} \left( \dsize \sum_{m<n}a_mx_m \right) = \dsize 
\sum_{m < n} a_mf(m)$ for any $a_m \in K^*$, so
$f^{\text{HOM}}$ is a functional from ${\Cal U}^*$ to $K^*$ and
$f^{\text{HOM}} \restriction {\Cal U}$ is functional from ${\Cal U}$ to $K$. 
\newline
We define support$(\dsize \sum_{m < \omega} a_mx_m) = \{m < \omega:a_m 
\ne 0_K\}$ for $\dsize \sum_{m < \omega} a_mx_m \in {\Cal U}^*$ and 
also content$(\dsize \sum_{m < \omega} a_mx_m) = \{a_m:m < \omega\}$. \newline
We define $P=P_{\Cal U}$ as the following simple forcing notion; so the 
set of members is defined by \scite{1.1}(a) where:
\medskip
\roster
\item "{$(a)$}"  $S = \{s:\text{for some } n,s \text{ is a function}$
\footnote{we could demand: if $n_1 < n_2 < n,s(n_\ell) = \dsize \sum_m a_mx_m
\in {\Cal U} \backslash \{0_{\Cal U}\}$; then max(support$(s(n_1)) <
\text{max(support}(s(n_2))\}$.  We define support
$( \dsize \sum_{m < \omega} a_mx_m) = \{m < \omega:a_m \ne 0_K\}$ for
$\dsize \sum_{m < \omega} a_mx_m \in {\Cal U}^*$ and also content
$(\dsize \sum_{m < \omega} a_m,x_m) = \{a_m:m < \omega\}$.
It appears more convenient to use the present definition.}
from $\{0,\dotsc,n-1\}$ to ${\Cal U}(\subseteq \omega)$ such that
$\langle s(\ell):\ell \in \text{ Dom}(s) \text{ and } s(\ell) \ne
0_{{\Cal U}^*} \rangle\}$ is linearly independent
\item "{$(b)$}"  $(s^1,\bold t) \le (s^2,\bold t^2)$ iff $s^1 \subseteq
s^2,\bold t^1 \subseteq \bold t^2$ and \newline
$[n \in \text{ Dom}(s^2) \backslash
\text{Dom}(s^1) \and t \in \bold t \and \text{ Rang}(t) \subseteq K
\Rightarrow t^{\text{HOM}}(s^2(n)) = 0]$.
\endroster
\medskip

\noindent
Here the following $f,g$ are O.K.: $f(n,m) = n + m + 1,g(s,\{t_\ell 
\restriction k:\ell < m\}) =$ first $x \in {\Cal U}^*$ under the order of the
natural numbers such that: $x$ has the form $\dsize \sum_{\ell < m} a_\ell 
x_\ell,x$ is not zero, and each $a_\ell$ is in the subfield $K_x$ of $K^*$ 
generated by $\dsize \bigcup_{\ell < k} \,\, 
\dsize \bigcup_{y \in \text{ Rang}(t_\ell)}$ content$(y) \cup 
\dsize \bigcup_{y \in \text{ Rang}(s)}$ content$(y)$
and $t^{\text{HOM}}_\ell(x) = 0$ for $\ell < m$.  This works as for any $n$
linear maps from a vector space over $K$ of dimension $n+m+1$ to $F$, the
intersection of their kernels has dimension at least $n+1$ and we can
restrict ourselves to the vector space over $K_x$ generated by
$\{x_\ell:\ell < n + m + 1\}$. \newline
4) Let $Q$ be defined by
\medskip
\roster
\item "{$(a)$}"  $S = \{s:s$ is a finite sequence, $s(\ell)$ is a finite
subset of $\omega,|s(\ell)| \le \ell + 1\}$
\smallskip
\noindent
\item "{$(b)$}"  $(s^1,\bold t^1) \le (s^2,\bold t^2)$ iff:
$s^1 \subseteq s^1,\bold t^1 \subseteq \bold t^2$ and \newline

$\qquad \qquad \qquad \qquad$ if $n \in \text{ Dom}(s^2) \backslash
\text{Dom}(s^1)$ and \newline

$\qquad \qquad \qquad \qquad s^2(n) \ne \emptyset$ then $(\forall t \in
\bold t^1)(t(n) \in s^2(n))$.
\endroster
\medskip

Note that this forcing suffices.  The value $\emptyset$ has little 
influence of our purpose.  Seemingly we can suppress it at the expense of a
slight burden on the iteration: making the $m$ in $(*)$ of \scite{1.14} be 
part of the condition; did not check.
\enddefinition   
\bigskip

\demo{\stag{1.3} Discussion}  1) We shall work with simple forcing as we 
cannot iterate Mathias forcing as this will make ${\frak p}$ large hence 
by Bell \cite{B} theorem MA($\sigma$-centered) holds, hence $\clubsuit$ 
fails so all this illustrates why we use the highly specialized ``simple 
forcing". \newline
2) Can we get $\clubsuit$ with otp$(C_\delta) > \omega$?  Well, we need to 
have $\theta_\xi > \omega$ and to allow to divide the set $\theta_\xi$ 
to finitely many convex parts, and this is done here. \newline
3) Can we get MA$_{\aleph_1}$[$\sigma$-centered]?  seems yes. \newline
4) We may use more strongly the $2^{\aleph_0} = \aleph_1$ in the end:
$k_\xi = 2 \Rightarrow$ the $u_\xi$ are submodel of some ${\frak B}$.
\newline
5) We can generalize replacing $\aleph_0$ by $\sigma = \sigma^{< \sigma}$,
so in \scite{1.4} below we demand $\kappa = \text{ cf}(\kappa) > \sigma,
\{1,2\} \subseteq \Theta \subseteq \sigma^+$, we use a strong version of 
$\sigma^+$-c.c. \newline
(e.g. as in \cite{Sh:80}).
\enddemo 
\bigskip

\centerline {$* \qquad * \qquad *$}
\bigskip

\noindent
The following defines a ``domain", ``base", ``carrier" of a condition in the
iterations defined in the Main Definition \scite{1.11} below.  For a more 
detailed explanation, see after the Definition.
\bigskip

\demo{\stag{1.4} Context}  1) $\kappa$ will be regular uncountable, can be
chosen as $\aleph_1$. \newline
2) $\Theta$ be a set of ordinals $\in (0,\kappa),1 \in \Theta,\Theta \ne
\{1\}$; now we shall need: $\Theta \subseteq \omega_1$ (in \scite{1.14}) and
$2 \in \Theta$ (or just $\Theta \cap [2,\omega) \ne 0$); in
\scite{1.19}(1) the case $\Theta = \{1,2,\omega\}$, really here
$V \models GCH,\kappa = \aleph_1,\mu = \aleph_2$ suffice for the main 
theorem (the other cases 
arise when we want to strengthen $\clubsuit$). \newline
3) cf$(\mu) \ge \kappa$.
\enddemo
\bigskip

\definition{\stag{1.5} Definition}  Let $\zeta < \kappa,\alpha$ an ordinal.  
An $(\alpha,\zeta,\Theta)$-creature $u$ is a sequence $\langle \zeta,\bar k,
\beta,\bar \theta,\bar \psi,\bar w,\bar F \rangle$ where (but we may
suppress $\Theta$, being constant):
\medskip
\roster
\item "{$(a)$}"  $\bar k = \langle k_\varepsilon:\varepsilon < \zeta \rangle$
where $k_\varepsilon \in \{1,2\}$
\smallskip
\noindent
\item "{$(b)$}"  $\bar \theta = \langle \theta_\varepsilon:\varepsilon <
\zeta \rangle$ satisfying: \newline
$[k_\varepsilon = 1 \Rightarrow \theta_\varepsilon = 1]$ and 
$[k_\varepsilon = 2 \Rightarrow \theta_\varepsilon > 1]$ \newline
$\theta_\varepsilon$ an ordinal $\in \Theta$
\smallskip
\noindent
\item "{$(c)$}"  $\bar \psi = \langle \psi_\varepsilon:\varepsilon < \zeta
\rangle$, where $\psi_\varepsilon < \theta_\varepsilon$
\smallskip
\noindent
\item "{$(d)$}"  $\bar w = \langle w_{\bar \varphi}:\bar \varphi
\in \Phi \rangle$ where:

$$
\Phi =: \bigcup \{\Phi_\xi:\xi \le \zeta \}
$$

$$
\align
\Phi_\xi =: \biggl\{ \bar \varphi:&\bar \varphi = \langle \varphi_\varepsilon:
\xi \le \varepsilon < \zeta \rangle,\varphi_\varepsilon < 
\theta_\varepsilon \\
  &\text{and dif} (\bar \psi,\bar \varphi) =: \{\varepsilon:\xi \le
\varepsilon < \zeta \text{ and } \psi_\varepsilon \ne \varphi_\varepsilon\} \\
  &\text{is finite} \biggr\}
\endalign
$$
\noindent
(we can think of $\bar \varphi \in \Phi_\xi$ as an alternate history of the
condition from time $\xi$ on)
\smallskip
\noindent
\item "{$(e)$}"  $\bar F = \langle F_{\bar \varphi_1,\bar \varphi_2}:
(\bar \varphi_1,\bar \varphi_2) \in \Psi \rangle$ where $\Psi = \dsize 
\bigcup_{\xi \le \zeta} \Psi_\xi$ where
$$
\Psi_\xi = \{(\bar \varphi^1,\bar \varphi^2):\bar \varphi^1 \in \Phi_\xi,
\bar \varphi^2 \in \Phi_\xi \text{ and } \bar \varphi^1 \restriction
[\xi + 1,\zeta) = \bar \varphi^2 \restriction [\xi + 1,\zeta)\}.
$$
\endroster
\medskip

\noindent
(note: if $k_\varepsilon = 1$ or $\xi = \zeta,\Psi_\xi$ has only pairs
of the form $(\bar \varphi,\bar \varphi)$ and if $\xi_1 \le \xi_2 \le \zeta,
\bar \varphi^2 \in \Phi_{\xi_2}$ then for some $\bar \varphi^1 \in
\Phi_{\xi_1}$ we have $\bar \varphi^2 = \bar \varphi^1 \restriction \xi_1$)
such that:
\medskip
\roster
\item "{$(f)$}"  if $\bar \varphi^1,\bar \varphi^2 \in \Phi$ and
$\bar \varphi^1 = \bar \varphi^2 \restriction [\xi,\zeta)$ \underbar{then}
$w_{{\bar \varphi}^2} \subseteq w_{{\bar \varphi}^1}$ and each
$w_{\bar \varphi}$ is a subset of $\alpha$ of cardinality $< \kappa$
\smallskip
\noindent
\item "{$(g)$}"  if $\bar \varphi \in \Phi_0$ and $\xi \le \zeta$ is a
limit ordinal \underbar{then} $w_{\bar \varphi \restriction \xi} =
\dsize \bigcup_{\varepsilon < \xi} w_{\bar \varphi \restriction
[\varepsilon,\zeta)}$
\smallskip
\noindent
\item "{$(h)$}"  if $\bar \varphi^1 \in \Phi_{\xi+1},\xi \le \zeta,
k_\xi = 2$ \underbar{then} \newline
$w_{\bar \varphi^1} = \bigcup\{ w_{{\bar \varphi}^2}:\bar \varphi^2 \in 
\Phi_\xi$ and $\bar \varphi^2 \restriction [\xi+1,\zeta) = \bar \varphi^1\}$
\smallskip
\noindent
\item "{$(i)$}"  for $(\bar \varphi^1,\bar \varphi^2) \in \Psi_\xi$ let
$\xi(\bar \varphi^1,\bar \varphi^2)$ be that (unique) $\xi$; for
$\bar \varphi \in \Phi_\xi$ let $\xi(\bar \varphi) = \xi$
\smallskip
\noindent
\item "{$(j)$}"  $F_{\bar \varphi^1,\bar \varphi^2}$ is a function,
Dom$(F_{\bar \varphi^1,\bar \varphi^2})$ is an initial segment of
$w_{\bar \varphi^2}$, Rang$(F_{\bar \varphi^1,\bar \varphi^2})$ is an
initial segment of $w_{\bar \varphi^1}$
\smallskip
\noindent
\item "{$(k)$}"  $F_{\bar \varphi^1,\bar \varphi^2}$ is a one to one and
order preserving function and $F_{\bar \varphi^2,\bar \varphi^1} =
F^{-1}_{\bar \varphi^1,\bar \varphi^2}$
\smallskip
\noindent
\item "{$(l)$}"  if $(\bar \varphi^1,\bar \varphi^2) \in \Psi,
(\bar \varphi^2,\bar \varphi^3) \in \Psi$ \underbar{then}
$F_{\bar \varphi^1,\bar \varphi^3} \supseteq F_{\bar \varphi^1,
\bar \varphi^2} \circ F_{\bar \varphi^2,\bar \varphi^3},
F_{\bar \varphi,\bar \varphi} =
\text{ id}_{w_{\bar \varphi}}$ for $\varphi \in \Phi_\xi$
\smallskip
\noindent
\item "{$(m)$}"  for $(\bar \varphi^1,\bar \varphi^2) \in \Psi_\xi$ let
$\bar \varphi^1 < \bar \varphi^2$ means $\varphi^1_\xi < \varphi^2_\xi$
\smallskip
\noindent
\item "{$(n)$}"  if $(\bar \varphi^1,\bar \varphi^2) \in \Psi_\xi$ and
$\bar \varphi^1 < \bar \varphi^2$ \underbar{then} $\gamma \in \text{ Dom}(F_
{\bar \varphi^1,\bar \varphi^2}) \Rightarrow F_{\bar \varphi^1
\bar \varphi^2}(\gamma) \le \gamma$
\smallskip
\noindent
\item "{$(o)$}"  if $\bar \varphi \in \Psi_{\xi+1},k_\xi = 2,\gamma <
\alpha$ and for at least two $\bar \varphi^1 \in \Phi_\xi$ we have \newline
$\bar \varphi^1 \restriction [\xi + 1,\zeta) = \bar \varphi \and
\gamma \in w_{\bar \varphi^1}$ \underbar{then}:
{\roster
\itemitem{ $(\alpha)$ }  for every $\bar \varphi^1 \in \Phi_\xi$ satisfying
$\bar \varphi^1 \restriction [\xi+1,\zeta) = \bar \varphi$ we have \newline
$\gamma \in w_{\bar \varphi^1}$
\itemitem{ $(\beta)$ }  if $\bar \varphi^1,\bar \varphi^2 \in \Phi_\xi$ and
$\bar \varphi^1 \restriction [\xi +1,\zeta) = \bar \varphi^2 \restriction
[\xi +1,\zeta) = \bar \varphi$ \underbar{then} \newline
$F_{\bar \varphi^1,\bar \varphi^2}(\gamma) = \gamma$ (so is well defined and
otp$(w_{\bar \varphi^1} \cap \gamma) = \text{ otp}(w_{\bar \varphi^2} 
\cap \gamma)$)
\endroster}
\item "{$(p)$}"  if $\bar \varphi \in \Phi_0$ then $w_{\bar \varphi}
= \emptyset$
\smallskip
\noindent
\item "{$(q)$}"  assume $\bar \varphi \in \Phi_{\xi+1},k_\xi = 2$;
there is $\bold e$, a convex \footnote{the convexity is natural, did not check
if necessary} equivalence relation on $\theta_\xi$ with finitely many
equivalence classes such that:
{\roster
\itemitem { $(\alpha)$ }  if $\bar \varphi^1,\bar \varphi^2 \in \Phi_\xi,
\bar \varphi^1 \restriction [\xi + 1,\zeta] = \bar \varphi^2 \restriction
[\xi + 1,\zeta)$ and $\varphi^1_\xi,\varphi^2_\xi$ are $\bold e$-equivalent
\underbar{then} $F_{\bar \varphi^1,\bar \varphi^2}$ is from 
$w_{\bar \varphi^1}$ onto $w_{\bar \varphi^2}$
\itemitem{ $(\beta)$ }  if $(\bar \varphi^1,\bar \varphi^2) \in \Phi_\xi,
(\bar \varphi^1,\bar \varphi^3) \in \Phi_\xi$ (hence
$\bar \varphi^1,\bar \varphi^3) \in \Phi_\xi)$ and $\bar \varphi^1 \bold e
\bar \varphi^2$ \underbar{then} $F_{\bar \varphi^1,\bar \varphi^3} =
F_{\bar \varphi^1,\bar \varphi^2} \circ F_{\bar \varphi^2,\bar \varphi^3}$
\newline
(note: if $(\forall \theta \in \Theta)(\theta \le \omega)$, which is the 
main case here, $\bold e$ can be replaced by $n < \omega$, with 
$\varphi^1_\xi \bold e \varphi^2_\xi$ being replaced by \newline 
$(n \le \varphi^1_\xi \and n \le \varphi^2_\xi) \vee (\varphi^1_\xi =
\varphi^2_\xi)$.)
\endroster}
\endroster
\enddefinition
\bigskip

\demo{\stag{1.6} Explanation}  The creature $u$ can be thought 
of as a part of the creation of a condition in Definition \scite{1.11}.  
For $\zeta = 0$ we have an
``atomic condition", we could make it e.g. speaking on one iterand but we have
made them empty.  If $\zeta = \xi +1,k_\zeta = 2$ we have $\theta_\xi$
conditions forming a $\triangle$-system of conditions of ``depth", ``length
of history" $\xi$ which we put together.  If $\zeta = \xi +1,k_\zeta = 1$ 
we extend a
condition of depth $\xi$ to satisfy a density requirement on the domain.  
If $\zeta$ is
limit, we take a limit of a ``nicely increasing sequence of conditions".
Now $\bar \psi$ is the ``history": if we have just put together $\theta_\xi$
condition of depth $\xi$ from a $\triangle$-system, it is a ``free" choice
which one lies on the
``main history line" and which are just ``alternate histories".  So
$\Phi_\xi$ is the set of ``possible" alternate histories from stage $\xi$
on $w_{\bar \varphi}$, and the domain of the condition which was the 
beginnings
of this history and $\bar F_{\bar \varphi^1,\bar \varphi^2}$ is an isomorphism
witnessing our having used a $\triangle$-system.  So for $\kappa > \aleph_1$,
we shall get even for e.g. $<_{\bold e}$ (see below) only strategic
completeness. \newline
Note: we tend to ignore the case $u^{[\bar \varphi^1]} = 
u^{[\bar \varphi^2]}$ (see \scite{1.1}(7) below), $\bar \varphi^1 \ne 
\bar \varphi^2$, anyhow all is preserved.
\enddemo
\bigskip

\definition{\stag{1.7} Definition}  1)  Let an $\alpha$-creature mean an 
$(\alpha,\xi)$-creature for some $\xi < \kappa$. \newline
2) Let CR$_{\alpha,\zeta}$ be the set of $(\alpha,\zeta)$-creatures,
CR$_\alpha = \dsize \bigcup_{\zeta < \kappa}$ CR$_{\alpha,\zeta}$.  If
$u \in \text{ CR}_\alpha$ let \newline
$u = \langle \zeta^u,\bar k^u,\bar \theta^u,
\bar \psi^u,\bar w^u,\bar F^u \rangle$ and $\zeta[u] = \zeta^u$, etc.
and let \newline
Dom$(u) = w^u = w[u] =: w^u_{<>}$; also $\Phi^u,\Phi^u_\xi,\Psi^u,
\Psi^u_\xi$ are defined accordingly. \newline
3) For $u_1,u_2 \in \text{ CR}_\alpha$ let $u_1 \le_e u_2$ \underbar{iff}
$\zeta[u_1] \le \zeta[u_2],\bar k^{u_1} = \bar k^{u_2} \restriction \zeta
[u_1]$, \newline
$\bar \theta^{u_1} = \bar \theta^{u_2} \restriction \zeta[u_2]$ and
for some $\bar \varphi^* \in \Phi^{u_2}_{\zeta[u_1]}$ we have:
\medskip
\roster
\item "{$(a)$}"  $\bar \psi^{u_1} \approx \bar \psi^{u_2} \restriction
\zeta[u_1]$ i.e. $\{\varepsilon:\varepsilon < \zeta[u_1]$ and
$\psi^{u_1}_\varepsilon \ne \psi^{u_2}_\varepsilon\}$ is finite
\item "{$(b)$}"  for $\bar \varphi \in \Phi^{u_1}$ we have
$w^{u_1}_{\bar \varphi} = w^{u_2}_{\bar \varphi \cup
\bar \varphi^*}$
\item "{$(c)$}"  for $(\bar \varphi^1,\bar \varphi^2) \in \Psi^{u_1}$ we
have $F^{u_1}_{\bar \varphi^1,\bar \varphi^2} = F^{u_2}_{\bar \varphi^1
\cup \bar \varphi^*,\bar \varphi^2 \cup \bar \varphi^*}$.
\endroster
\medskip

\noindent
4) For $u_1,u_2 \in \text{ CR}_\alpha$ let $u_1 \le_{\text{de}} u_2$
\underbar{iff} $u_1 \le_e u_2$ and $\bar \psi^{u_1} = \bar \psi^{u_2}
\restriction \zeta[u_1]$ so in (3), $\bar \varphi^* = \bar \psi^{u_2}
\restriction [\zeta[u_1],\zeta[u_2])$. \newline
5) If $u \in \text{ CR}_\alpha$ and $\beta \le \alpha$ we define
$u \restriction \beta$ as $\langle \zeta^u,\bar k^u,\bar \theta^u,
\bar \psi^u,\bar w',\bar F' \rangle$, where \newline
$\bar w' = \langle \bar w^u_{\bar \varphi} \cap \beta:\bar \varphi \in
\Phi^u \rangle$ and $\bar F' = \langle F^u_{\bar \varphi^1,\bar \varphi^2} 
\cap (\beta \times
\beta):(\bar \varphi^1,\bar \varphi^2) \in \Psi^u \rangle$ (see Claim
\scite{1.9}(1). \newline
6) For $u_1,u_2 \in \text{ CR}_\alpha$ we say $u_1 \approx u_2$ iff
$u_1 \le_e u_2 \le_e u_1$. \newline
7) If $u \in \text{ CR}_{\alpha,\zeta},\xi < \zeta$ and $\bar \varphi \in
\Phi^u_\xi$ \underbar{then} let $v = u^{[\bar \varphi]} \in \text{ CR}
_{\alpha,\xi}$ be defined by:

$$
\zeta^v = \xi
$$

$$
\bar k^v = \bar k^u \restriction \xi
$$

$$
\bar \theta^v = \bar \theta^u \restriction \xi
$$

$$
\bar \psi = \bar \psi^u \restriction \xi
$$

$$
w^v_{\bar \varphi^1} = w^u_{\bar \varphi^1 \cup \bar \varphi}
$$

$$
F^v_{\bar \varphi^1,\bar \varphi^2} = F^u_{\bar \varphi^1 \cup
\bar \varphi,\bar \varphi^2 \cup \bar \varphi}.
$$
\medskip

\noindent
8) If $u \in \text{ CR}_{\alpha,\zeta}$ and $\xi \le \zeta^u$ \underbar{then}
$v = u^{[\xi]}$ is defined as $u^{[\bar \psi^u \restriction [\xi,\zeta)]}$.
\newline
9) If $u \in \text{ CR}_{\alpha,\zeta},\zeta^u = \xi + 1,k^u_\xi = 1$ 
\underbar{then} we let $u^{[*]} = u^{[\xi]}$. \newline
10)  If $u \in \text{ CR}_\alpha,\zeta^u = \xi +1,k^u_\xi = 2$
\underbar{then} let $\theta^{(u)} = \theta^u_\xi$ and $u^{(i)} =
u^{[\{ \langle \xi,i \rangle\}]}$ for $i < \theta^{(u)}$ and
$F^u_{(i,j)} = F^u_{\{ \langle \xi,i \rangle\},\{\langle \xi,j \rangle\}}$.
If $\xi < \zeta^u,k^u_\xi = 2,i < \theta^u_\xi$ \underbar{then} let \newline
$u^{[\xi],(i)} = (u^{{[\xi]})^{(i)}}$. \newline
11) For $u,v \in \text{ CR}_\alpha$ we say that $F$ maps $u$ to $v$
\underbar{if}: $F \restriction w^u$ is a one to one order preserving
function form $w^u$ onto $w^v,\zeta^u = \zeta^v,\bar k^u =
\bar k^v,\bar \theta^u = \bar \theta^v,\bar \psi^u = \bar \psi^v$, \newline
$w^v_{\bar \varphi} = F'' w^u_{\bar \varphi},F^v_{(\bar \varphi^1,
\bar \varphi^2)} = F \circ F^u_{(\bar \varphi^1,\bar \varphi^2)} \circ
F^{-1}$ and $F^u_{(\bar \varphi^1,\bar \varphi^2)} = F^{-1} \circ
F^v_{(\bar \varphi^1,\bar \varphi^2)} \circ F$.  We write this
$F(u) = v$.  \newline
12) For $u \in \text{ CR}_{\alpha,\zeta},\zeta \le \zeta^u,\{\bar \varphi^1,
\bar \varphi^2\} \subseteq \Phi^u_\zeta$, we define $F^u_{\bar \varphi^1,
\bar \varphi^2}$, an order preserving function from an initial segment of 
Dom$(u^{[\bar \varphi^2]})$ onto some initial segment of 
Dom$(u^{[\bar \varphi^1]})$, by induction on $\zeta^u$:
\medskip
\roster
\item "{$(a)$}"  if $\zeta = \zeta^u$ then $\bar \varphi^1 = \bar \varphi^2 =
\langle \rangle$ and we let $F^u_{\bar \varphi^1,\bar \varphi^2} =
\text{ id}_{\text{Dom}(u)}$
\medskip
\noindent
\item "{$(b)$}"  if $\zeta < \zeta^u = \xi + 1,k^u_\xi = 1$ we let
$F^u_{\bar \varphi^1,\bar \varphi^2} = F^{u^{[*]}}_{\bar \varphi^1 
\restriction [\zeta,\xi),\bar \varphi^2 \restriction [\zeta,\xi)}$
\medskip
\noindent
\item "{$(c)$}"  if $\zeta < \zeta^u = \xi +1,k^u_\xi = 2$ we let
$F^u_{\bar \varphi^1,\bar \varphi^2} = F^{u^{[\varphi^1_\xi]}}
_{\bar \varphi^1 \restriction [\zeta,\xi),\bar \varphi^2 \restriction 
[\zeta,\xi)} \circ F^u_{\varphi^1_\xi,\varphi^2_\xi}$
\medskip
\noindent
\item "{$(d)$}"  if $\zeta < \zeta^u,\zeta^u$ a limit ordinal let
$F^u_{\bar \varphi^1,\bar \varphi^2} = F^{u^{[\xi]}}_{\bar \varphi^1
\restriction [\zeta,\xi),\bar \varphi^2 \restriction [z\eta,\xi)}$ for every
$\xi < \zeta^u$ that is large enough such that
$\bar \varphi^1 \restriction [\xi,\zeta^u] = \bar \varphi^2 \restriction
[\xi,\zeta^u] = \bar \psi^u \restriction [\xi,\zeta^u)$.
\endroster
\enddefinition
\bigskip

\definition{\stag{1.8} Definition}  1) An $\alpha$-sequence is 
$\bar a = \langle a_\beta:\beta < \alpha \rangle$ such that $a_\beta 
\subseteq \beta$.  We say
$\bar a$ is a $(\mu,\alpha)$-sequence if in addition $|\alpha_\beta| < \mu$ 
for $\beta < \alpha$. \newline
2) For an $\alpha$-sequence $\bar a$ let CR$_{\bar a,\zeta}$ be the set of
$u \in \text{ CR}_{\alpha,\zeta}$ such that:
\medskip
\roster
\item "{$(*)$}"  if $\xi < \zeta^u$ and $k_\xi = 2$ and $\theta^u_\xi \ge 
\omega$ and $(\bar \varphi^1,\bar \varphi^2) \in \Psi_\xi$ and $\varphi^1_\xi
< \varphi^2_\xi$ \newline
\underbar{then} 
$w_{\bar \varphi^2} \cap \cup \{a_\beta:\beta \in w_{\bar \varphi^1} \cap
w_{\bar \varphi^2}\} \subseteq w_{\bar \varphi^1}$.
\endroster
\medskip

\noindent
In this case, we say $u$ is an $(\bar a,\zeta)$-creature and we say $u$ is an
$\bar a$-creature.
\enddefinition
\bigskip

\proclaim{\stag{1.9} Claim}  1) If $u$ is an 
$(\alpha,\zeta)$-creature, $\beta \le \alpha$ \underbar{then} 
$u \restriction \beta$ (see Definition \scite{1.7}(5)) is a 
$(\beta,\zeta)$-creature,
(and $\zeta[u \restriction \beta] = \zeta[u]$ of course).  If $\bar a$ is an
$\alpha$-sequence, $u$ is an $(\alpha,\zeta)$-creature, $\beta \le \alpha$
then $u \restriction \beta$ is an $(\bar a \restriction \beta,
\zeta)$-creature.  If $\gamma \le \beta \le \alpha$ and $u$ is an
$(\alpha,\zeta)$-creature then $u \restriction \gamma = (u \restriction \beta)
\restriction \gamma$.  If $u$ is a $(\beta,\zeta)$-creature and
$\beta \le \alpha$ then ($u$ is an $(\alpha,\zeta)$-creature and)
$u \restriction \beta = u$.  \newline
2)  If $u$ is an $\alpha$-creature \underbar{then} for a unique $\zeta =
\zeta^u < \kappa,u$ is an $(\alpha,\zeta)$-creature. \newline
3) $\le_e$ is a partial order on CR$_\alpha,[u_1 \le_e u_2 \Rightarrow
\zeta(u_1) \le \zeta(u_2)]$ and $[u_1 \le_e u_2 \and \zeta[u_1] = \zeta[u_2]
\Leftrightarrow u_1 \approx u_2]$ and $[u_1 \le_e u_2 \Leftrightarrow (\exists
\bar \varphi \in \Phi^{u_2})(u_1 \approx u_2^{[\bar \varphi]})]$ and for
an $\alpha$-sequence $\bar a$ we have: $u_1 \le_e u_2 \and u_2 \in
\text{ CR}_{\bar a} \Rightarrow u_1 \in \text{ CR}_{\bar a}$. \newline
4) $\le_{\text{de}}$ is a partial order on CR$_\alpha,[u_1 \le_{\text{de}}
u_2 \Rightarrow u_1 \le_e u_2]$ and $[u_1 \le_{\text{de}} u_2 \and
\zeta[u_1] = \zeta[u_2] \Rightarrow u_1 = u_2]$ and $[u_1 \le_{\text{de}}
u_2 \Leftrightarrow (\exists \xi \le \zeta^{u_2})[u_1 = u_2^{[\xi]}]$.
\newline
5) $u^{[\xi]}$ is the unique $v$ such that $\zeta^v = \xi \and v 
\le_{\text{de}} u$ (defined iff $\xi \le \zeta^u$), and for $\xi \le
\varepsilon \le \zeta^u,(u^{[\varepsilon]})^{[\xi]} = u^{[\xi]}$.  Also for
$u \in \text{ CR}_{\alpha,\zeta}$, if $\xi \le \zeta$ and $\bar \varphi \in
\Phi^u_\xi$ \underbar{then} $u^{[\bar \varphi]} \in \text{ CR}_{\alpha,\xi}$
and $u^{[\bar \varphi]} \le_e u$. \newline
6) On CR$_\alpha,u \approx u'$ is an equivalence relation (see
Definition \scite{1.7}(6)). \newline
7) If $u_1 \le_e u_2$ \underbar{then} for some $u'_2 \approx u_2$ we have
$u_1 \le_{\text{de}} u'_2$. \newline
8) If $u^1 \le_e u^2$ and $u_1 \le_e u^1 \restriction \beta$ \underbar{then}
for some $u_2 \le_e u^2$ we have $u_1 = u_2 \restriction \beta$. \newline
9) If $\beta < \alpha$ \underbar{then} CR$_{\beta,\zeta} \subseteq
\text{ CR}_{\alpha,\zeta}$ and for any $\alpha$-sequence $\bar a$,
CR$_{\bar a \restriction \beta,\zeta} \subseteq \text{ CR}_{\bar a,\zeta}$.
\newline
10) If $u \in \text{ CR}_{\alpha,\zeta},\zeta = \xi +1,k^u_\xi = 1$
\underbar{then} $u^{(*)} = u^{[\xi]} = u^{[\langle \xi,0 \rangle]}$.
\newline
11) If $u \in \text{ CR}_{\alpha,\zeta},\zeta = \xi +1,k^u_\xi = 2$
\underbar{then} $\{v^{[\xi]}/\approx:v \approx u\} = \{u^{(i)}/\approx:i <
\theta^u_\xi\}$. \newline
12) If $\delta < \kappa$ is a limit ordinal, $\langle u_i:\alpha < \delta 
\rangle$ is a $\le_{\text{de}}$-increasing sequence in CR$_\alpha$ 
\underbar{then} for one and only one $u \in \text{ CR}_\alpha$ we have:
$[i < \delta \Rightarrow u_i
\le_{\text{de}} u]$ and $\zeta^u = \dsize \bigcup_{i < \delta} \zeta^{u_i}$. 
\newline
13) For $u \in CR_\alpha,\xi < \zeta^u$ and $\bar \varphi^1,\bar \varphi^2,
\bar \varphi^3 \in \Phi^u_\xi$ we have:
\medskip
\roster
%
\item "{$(i)$}"  $F^u_{\bar \varphi^1,\bar \varphi^2}$ (defined in
\scite{1.7}(12)) is really a one to one order preserving function from some
initial segment of $w^u_{\bar \varphi^2}$ onto some initial segment of
$w^u_{\bar \varphi^1}$
\item "{$(ii)$}"  $\gamma \in \text{ Dom}(F^u_{\bar \varphi^1,\bar \varphi^2})
\cap \text{ Rang}(F^u_{\bar \varphi^1,\bar \varphi^2}) \Rightarrow
F^u_{\bar \varphi^1,\bar \varphi^2}(\gamma) = \gamma$
\item "{$(iii)$}"  $F^u_{\bar \varphi^2,\bar \varphi^1} = 
(F^u_{\bar \varphi^1,\bar \varphi^2})^{-1}$
\item "{$(iv)$}"  $F^u_{\bar \varphi^1,\bar \varphi^3} \supseteq
F^u_{\bar \varphi^1,\bar \varphi^2} \circ F^u_{\bar \varphi^2,\bar \varphi^3}$
\item "{$(v)$}"  Dom$(F^u_{\bar \varphi^1,\bar \varphi^2}) = 
w^u_{\bar \varphi^2} \Leftrightarrow \text{ Rang}(F^u_{\bar \varphi^1,
\bar \varphi^2}) = w^u_{\bar \varphi^1}$
\item "{$(vi)$}"  if in (v) the equalities hold \underbar{then} \newline
$F^u_{\bar \varphi^1,\bar \varphi^2}(u^{[\bar \varphi^2]}) =
u^{[\bar \varphi^1]}$
\item "{$(vii)$}"  if in (v) the equalities fail and \newline
$\gamma_2 = \text{ Min}(w^u_{\bar \varphi^2} \backslash
\text{ Dom}(F^u_{\bar \varphi^1,\bar \varphi^2}))$ and \newline
$\gamma_1 = \text{ Min}(w^u_{\bar \varphi^1} \backslash \text{ Dom}
(F^u_{\bar \varphi^2,\bar \varphi^1}))$ \underbar{then}
$F^u_{\bar \varphi^1,\bar \varphi^2}(u^{[\bar \varphi^2]} \restriction
\gamma_2) = u^{[\bar \varphi^1]} \restriction \gamma_1$.
\endroster
\medskip

\noindent
14) Assume $u \in \text{ CR}_{\alpha,\zeta},v \le_e u,u_1 = 
u^{[\bar \varphi]},\bar \varphi \in \Phi^u_{\xi +1},k^u_\xi = 2,\theta =
\theta^u_\xi$, \newline
$v_i = u^{(i)}_1(= u^{[\langle i \rangle \char 94 \bar \varphi]}$,
see Definition \scite{1.7}(10)). \newline
\underbar{Then} there is a convex equivalence relation $\bold e$ on 
$\theta$ with finitely many equivalence classes such that:
\medskip
\roster
\item "{$(*)$}"  $i \bold e j$ implies:
{\roster
\itemitem{ $(a)$ }  $F^{u_1}_{(i,j)}$ is an order preserving map from
Dom$(v_i)$ onto Dom$(v_j)$
\itemitem{ $(b)$ }  $F^{u_1}_{(i,j)}$ maps Dom$(v) \cap \text{ Dom}(v_i)$ onto
Dom$(v) \cap \text{ Dom}(v_j)$.
\endroster}
\endroster
\medskip

\noindent
15)  If $\gamma \in \text{ Dom}(u),u \in CR_\alpha$, \underbar{then} 
for a unique $\xi < \zeta^u$ we have: $k_\xi = 1$, and $\Phi^{u,\gamma} 
=: \{\bar \varphi \in \Phi^u_{\xi +1}:\gamma \in w^u_{\bar \varphi}\}$ 
is non-empty but
$[\bar \varphi \in \dsize \bigcup_{\varepsilon \le \xi} \Phi_\varepsilon
\Rightarrow \gamma \notin w^u_{\bar \varphi}]$ and
$[\bar \varphi \in \dsize \bigcup_{\varepsilon > \xi} \Phi^u_\varepsilon
\Rightarrow (\exists \bar \varphi^1)(\bar \varphi^1 \in \Phi^{u,\gamma} \and
\bar \varphi = \bar \varphi^1 \restriction \xi)$ and $|a_\gamma| < \kappa
\and \bar \varphi^1 \in \Phi^{u,\gamma} \and \bar \varphi^2 \in
\Phi^{u,\gamma} \Rightarrow w^u_{\bar \varphi^1} \cap a_\gamma =
w^u_{\bar \varphi^2} \cap a_\gamma$. \newline
16)  Assume $\bar \varphi \in \Phi^u_{\xi +1},k_\xi=2$ and
$\bar \varphi^1,\bar \varphi^2 \in \Phi^u_\xi$, \underbar{then} \newline 
$w_{\bar \varphi^1} \cap w_{\bar \varphi^2} = \{\gamma:
F^u_{\bar \varphi^1,\bar \varphi^2}(\gamma) = \gamma\} = \text{ Dom}
(F^u_{\bar \varphi^1,\bar \varphi^2}) \cap \text{ Dom}(F^u_{\bar \varphi^2,
\bar \varphi^1})$. 
\endproclaim
\bigskip

\demo{Proof}  Straightforward, but we elaborate 1), 8) and 14).
\enddemo
\bigskip

\noindent
1) Reading Definition 1.5, the least trivial clause to check is (q), so
assume \newline
$\xi < \zeta^{u \restriction \beta},\bar \varphi \in \Phi
^{u \restriction \beta}_{\xi + 1}$ and $k^{u \restriction \beta}_\xi =2$.
But $\zeta^{u \restriction \beta} = \zeta^u,\Phi^{u \restriction \beta}
_{\xi +1} = \Phi^u_{\xi +1}$ and \newline
$k^{u \restriction \beta}_\xi = k^u_\xi$ and
$\theta^{u \restriction \beta}_\xi = \theta^u_\xi$, so as $u \in CR_\alpha$,
there is a convex equivalence relation $\bold e^1$ on $\theta^u_\xi$ as
required on $\bold e$ in clause (q) of Definition \scite{1.5} for $u$.  We
now define another equivalence relation $\bold e^2$ on $\theta^u_\xi$;
letting

$$
\align
i \bold e^2 j \text{ \underbar{iff} } &i,j < \theta^{u,i \bold e^1 j}
\text{ and} \\
 &(\forall \gamma \in w_{\bar \varphi \cup \{ \langle \xi,i \rangle\}})
[\gamma < \beta \leftrightarrow F^u_{\bar \varphi \cup \{ \langle \xi,j 
\rangle\},\bar \varphi \cup \{ \langle \xi,i \rangle\}}(\gamma) < \beta].
\endalign
$$
\medskip

\noindent
Clearly $\bold e^2$ is an equivalence relation on $\theta^u_\xi$. \newline
By clause (n) of Definition \scite{1.5} clearly $\bold e^2$ is convex.

Now if $\bold e^2$ has infinitely many equivalence classes, let $j$ be
minimal such that $\{i/\bold e^2:i \bold e^1 j\}$ is infinite.  Let $i_n$ be
the first element in the $n$-th $\bold e^2$-equivalence class 
$\subseteq j/\bold e^1$
(or just $i_n < i_{n+1},\neg i_n \bold e^2 i_{n+1},i_n \bold e^1 j$) and let
$\gamma_n \in w_{\bar \varphi \cup \{ \langle \xi,i_n \rangle\}}$ witness
$\neg i_n \bold e^2 i_{n+1}$; i.e. $\gamma'_n =: F^u_{\bar \varphi \cup \{ 
\langle \xi,i_{n+1} \rangle\},\bar \varphi \cup \{ \langle \xi,i_n \rangle\}}
(\gamma_n) < \beta \nleftrightarrow \gamma_n < \beta$.

Without loss of generality $\gamma_n$ is minimal under this condition, so by
clause (n) of Definition \scite{1.5} as $i_n < i_{n+1}$ we have
$\gamma_n < \beta \le \gamma'_n$, hence necessarily $\gamma_{n+1} < \gamma'_n$
hence $F^u_{\bar \varphi \cup \{ \langle \xi,i_n \rangle\},\bar \varphi \cup 
\{ \langle \xi,i_{n+1} \rangle\}}(\gamma_{n+1}) < \gamma_n$ (remember that
by clause $(q)(\alpha)$ of Definition \scite{1.5},
$F^u_{\bar \varphi \cup \{ \langle \xi,i_n \rangle\},\bar \varphi \cup 
\{ \langle \xi,i_{n+1} \rangle\}}$ is from $w^u_{\bar \varphi \cup \{ 
\langle \xi,i_{n+1} \rangle\}}$ onto $w^u_{\bar \varphi \cup \{ \langle 
\xi,i_n \rangle\}}$).  Hence $\bigl< F^u_{\bar \varphi \cup 
\{ \langle \xi,i_0 \rangle\},\bar \varphi \cup \{ \langle \xi,i_n \rangle\}}
(\gamma_n):n < \omega \bigr>$ is a strictly decreasing sequence of ordinals;
contradiction. \newline
8) Let $\bar \varphi \in \Phi^{u^2}$ be such that $u^1 = 
(u^2)^{[\bar \varphi]}$ and let $\bar \varphi^1 \in \Phi^{u^1 \restriction
\beta}$ be such that $u_1 = (u^1 \restriction \beta)^{[\bar \varphi^1]}$
(exists by the definition of $\le_e$ and the assumptions). \newline
Clearly Dom$(\bar \varphi) = [\zeta^{u^1},\zeta^{u^2})$ and
Dom$(\bar \varphi^1) = [\zeta^{u_1},\zeta^{u^1 \restriction \beta})$.
However, $\zeta^{u^1 \restriction \beta} = \zeta^{u^1},\Phi^{u^1 \restriction
\beta} = \Phi^{u^1}$ hence $\bar \varphi^1 \in \Phi^{u^1}$, so easily
$\bar \varphi^1 \cup \bar \varphi \in \Phi^{u^2}$ and Dom$(\bar \varphi^1 \cup
\bar \varphi) = [\zeta^{u_1},\zeta^{u^2})$, lastly let $u_1 = (u^2)
^{[\bar \varphi^1 \cup \bar \varphi]}$; check.
\bigskip

\noindent
14)  We prove it by induction of $\zeta^u$.  Let $\bold e^1$ be a convex
equivalence relation on $\theta$ with finitely many equivalence classes which
satisfies clause $(a)$ of $(*)$ (exists by clause $(q)$ of Definition
\scite{1.5}).  We will refine it to satisfy clause $(b)$ too. \newline
If $v = u$ the conclusion is immediate: use $\bold e = \bold e^1$, so we
shall assume $v \ne u$.
\bigskip

\noindent
\underbar{Case 1}:  $\zeta^u = 0$.

Nothing to prove.
\bigskip

\noindent
\underbar{Case 2}:  $\zeta^u$ successor.
\bigskip

\noindent
\underbar{Subcase 2A}:  $\zeta[u_1] = \zeta^u$, so $\xi = \zeta^u - 1$.
\newline
Let $\bar \varphi^1$ be such that $v = u^{[\bar \varphi^1]},
\bar \varphi^1 \in \Phi^u$, note (as $v \ne u$) that necessarily 
Dom$(\bar \varphi^1) \ne \emptyset$, so $\xi \in \text{ Dom}
(\bar \varphi^1)$ and let $i = \varphi^1_\xi \,(< \theta)$ so 
$v \le_e u^{(i)}$.
Define $\bold e^2$, an equivalence relation on $\theta$ with the
classes $\{j:j < i\},\{i\},\{j:i < j < \theta\}$ (omitting any occurence of
the empty set).  Let $\bold e = \bold e^1 \cap \bold e^2$.  By clauses
$(n),(o)$, of Definition \scite{1.5} we are done.
\bigskip

\noindent
\underbar{Subcase 2B}:  $\zeta(u_1) < \zeta^u$.

If $k^u_{\zeta[u]-1} = 1$, necessarily $v \le_e u^{[*]}$ (see Definition
\scite{1.7}(9)), and $u_1 \le u^{[*]}$ and $\zeta[u^{[*]}] < \zeta^u$ so
we can use the induction hypothesis.  So assume $k_{\zeta[u]-1} = 2$, hence
for some $i^0,i^1 < \theta^u_{\zeta[u]-1}$ we have $v \le_e u^{(i^0)},
u_1 \le_e u^{(i)}$.  If
$i^0 = i^1$ as $\zeta[u^{(i^0)}] = \zeta^u -1$ we can use the induction
hypothesis; so assume $i_0 \ne i_1$.  Let $v = u^{[\bar \varphi^1]}$ so
$\varphi^1_{\zeta[u]-1} = i_0$, let
$\bar \varphi^2 = \langle \varphi^2_\zeta:\zeta^v \le \zeta < \zeta^u \rangle$
be: $\varphi^2_\zeta$ is $\varphi^1_\zeta$ if $\zeta \in [\zeta^v,\zeta^u
-1)$, and $i^1$ if $\zeta = \zeta^u -1$.  Now apply the induction hypothesis
with $u^{(i_1)},u_1,u^{[\bar \varphi^2]}$ here standing for $u,u_1,v$ there
and get an equivalence relation $\bold e$.  Clearly $F^u_{(i^1,i^0)}$
maps $v$ to $u^{[\bar \varphi^2]}$.  Is $\bold e$ as required?  So
assume $j^1,j^2 < \theta,j_1 \bold e j_2$ and $\gamma \in \text{ Dom}(v) \cap
\text{ Dom}(u^{(j^1)}_1)$; then necessarily $\gamma \in \text{ Dom}(u^{(i^0)})
\cap \text{ Dom}(u^{(i^1)})$ hence $F^u_{(i^1,i^0)}(\gamma) = \gamma$ (see
Definition \scite{1.5} clause $(o)$) so $\gamma \in 
\text{ Dom}(u^{[\bar \varphi^2]})$ hence $\gamma \in \text{ Dom}(
u^{[\bar \varphi^2]}) \cap \text{ Dom}(u^{(j^1)}_1)$ hence by the choice
of $\bold e$ we have $\gamma \in \text{ Dom}(u^{[\bar \varphi^2]}) \cap 
\text{ Dom}(u^{(j^2)}_1)$ hence $\gamma \in \text{ Dom}(u^{(j^2)}_1)$ 
hence $\gamma \in
\text{ Dom}(v) \cap \text{ Dom}(u^{(j^2)}_1)$.  By the symmetry between
$j^1$ and $j^2$ this suffices.
\bigskip

\noindent
\underbar{Case C}:  $\zeta^u$ is a limit ordinal.

We are assuming $\zeta^v < \zeta^u$ and as $\zeta^{u_1}$ is a successor
ordinal $\le \zeta^u$ it is $< \zeta^u$ so for every $\xi < \zeta^u$ be 
large enough,
$v,u_1 \le_e u^{[\xi]}$, and we can use the induction hypothesis. \newline
${}$ \hfill$\square_{\scite{1.9}}$
\bigskip

\centerline {$* \qquad * \qquad *$}
\bigskip

\definition{\stag{1.10} Definition/Claim}  1) If $\beta \le \alpha,
u_1 \in \text{ CR}_\alpha,u_0 =: u_1 \restriction \beta \le_{\text{de}}
u_2 \in \text{ CR}_\beta$ we defined $u = u_1 \otimes u_2 \in 
\text{ CR}_\alpha$ as follows.  First
\medskip
\roster
\item "{$(a)$}"   $\zeta^u = \zeta^{u_2},\bar k^u = \bar k^{u_2},
\bar \theta^u = \bar \theta^{u_2}$ and $\bar \psi^u = \bar \psi^{u_2}$.
\endroster
\noindent
Second
\roster
\item "{$(b)$}"  $w^u_{\bar \varphi}$ is defined as follows:
{\roster
\smallskip
\noindent
\itemitem { $(i)$ }    if $\xi(\bar \varphi) \le \zeta^{u_1}$ and
$\bar \varphi \restriction [\zeta^{u_1},\zeta^{u_2}) = \psi^{u_2} 
\restriction [\zeta^{u_1},\zeta^{u_2})$ then \newline
$w^u_{\bar \varphi} = w^{u_1}_{\bar \varphi \restriction 
((\text{Dom }\bar \varphi) \cap [0,\zeta^{u_1}))}$ 
\smallskip
\noindent
\itemitem { $(ii)$ }  if $\xi \in (\zeta^{u_1},\zeta^{u_2}]$ and 
$\bar \varphi = \psi^{u_2} \restriction [\xi,\zeta^{u_2})$ then \newline
$w^u_{\bar \varphi} = w^{u_2}_{\bar \varphi} \cup w^{u_1}$
\smallskip
\noindent
\itemitem { $(iii)$ }  if both cases do not apply then $w^u_{\bar \varphi} = 
w^{u_2}_{\bar \varphi}$.
\endroster}
\endroster
\noindent
Lastly
\roster
\item "{$(c)$}"  $F^u_{\bar \varphi^1,\bar \varphi^2}$ are defined naturally.
\endroster
\medskip

\noindent
2)  In part (1), $u$ actually belongs to CR$_\alpha$ and
$u \in \text{ CR}_\alpha,u \restriction \beta = u_2$ and
$u_1 \le_{\text{de}} u$. \newline
3)  If $u_1 = u_2 \restriction \beta$ and $\beta' = \cup\{ \gamma +1:\gamma
< \beta \text{ and } \gamma \in \text{ Dom}(u_2)$ (equivalently $\gamma \in
\text{ Dom}(u_1)\}$ then $u_1 = u_2 \restriction \beta'$. \newline
4)  Let $u_1 \le_{rd} u_2$ means: for some $\beta$ and $u$ we have $u_1 = u
\restriction \beta,u \le_e u_2$.  We call such $(\beta,u)$ a witnessing 
pair to $u_1 \le_{rd} u_2$, and if $\beta$ is minimal (as in (3)) we call 
such $(\beta,u)$ a good witnessing pair in $u_1 \le_{rd} u_2$. \newline
5) $\le_{rd}$ is partial order (on creatures).  In fact, if $u_0 \le_{rd}
u_1 \le_{rd} u_2$ and $(\beta_\ell,u'_\ell)$ is the witnessing pair for
$u_\ell \le_{rd} u_{\ell +1}$ and $u''_0 \le_e u'_1$ is such that
$u_0 = u''_0 \restriction \beta_0$ (exist by \scite{1.9}(8)) \underbar{then}
$(\beta_0,u''_0)$ is the witnessing pair for $u_0 \le_{rd} u_2$.
\enddefinition
\bigskip

\demo{Proof}  Straightforward.
\enddemo
\bigskip

\centerline{$* \qquad * \qquad *$}
\bigskip

\noindent
We turn to the iteration
\proclaim{\stag{1.11} Definition/Lemma}   We define and prove the following by
induction on $\alpha$ for any $\alpha$-sequence $\bar a$:
\medskip
\roster
\item "{$(A)$}"  We define when $\bar Q = \langle P_\beta,
{\underset\sim {}\to Q_\beta}:\beta < \alpha \rangle$ is (an iteration) from
\footnote{Remark: Note that e.g. $P_\beta * {\underset\sim
{}\to Q_\beta} \le P_{\beta +1}$ is not included, in fact, this does not
really interest us.  What we want is condition (b) of Conclusion 
\scite{1.18}(1) + \scite{1.19}(3).} $K^{\bar a}$ 
\item "{$(B)$}" for $\bar Q \in K^{\bar a}$ we define the set of elements of
\newline
$P_\alpha = \text{ Lim}(\bar Q) = \text{ Lim}^\alpha(\bar Q)
= \text{ Lim}^{\alpha,\bar a}(\bar Q)$
\item "{$(C)$}"  for $\bold p \in P_\alpha$ and $\beta \le \alpha$ we define
$\bold p \restriction \beta \in P_\beta,u[\bold p]$ and also $\bold p^{[u]}
\in P_\alpha$ for $u \le_e u[\bold p]$ (and show $\bold p \restriction \beta$ 
is well defined and in $P_\beta,u[\bold p]$ is well defined and in CR$_\alpha$
and $\bold p^{[u]}$ is well defined and in $P_\alpha$).  We also define when
$F(\bold p) = \bold q$,
\item "{$(D)$}"  we define the partial order $\le^{P_\alpha} = <^{P_\alpha}
_{nr}$ on $P_\alpha$ as well as $\le^{P_\alpha}_{pr},\le^{P_\alpha}_{apr},
<^{P_\alpha}_{vpr}$ \newline
(we omit the superscript $P_\alpha$ when clear)
\item "{$(E)$}"  we prove: the relations $\le^{P_\alpha},\le^{P_\alpha}_{pr},
\le^{P_\alpha}_{vpr}$ and $\le^{P_\alpha}_{apr}$ are partial orders on
$P_\alpha$ and if $\gamma \le \beta \le \alpha$ \underbar{then} $P_\gamma
\subseteq P_\beta \subseteq \text{ Lim}^\alpha(\bar Q)$ and $\le^{P_\beta}
\restriction P_\gamma = \le^{P_\gamma}$; similarly for $\le_{pr},\le_{apr}$
and $\le_{vpr}$; also for $\bold p \in P_\alpha,\bold p \restriction \gamma
= (\bold p \restriction \beta) \restriction \gamma$, \newline
$\bold p \restriction \gamma \le^{P_\alpha}_x \bold p \restriction 
\beta,[\bold p \le^{P_\beta}_x \bold q
\Rightarrow \bold p \restriction \gamma \le^{P_\gamma}_x \bold q \restriction
\gamma]$ (for $x \in \{nr,pr,vpr,apr\}$) and 
$\bold p \le^{P_\alpha}_x \bold q \Rightarrow \bold p \le^{P_\alpha}_y 
\bold q$ if \newline
$(x,y) \in \{(vpr,nr)\},(pr,nr),(apr,nr),(vpr,pr)\}$
\item "{$(F)$}"  for $\bold p \in P_\alpha,\beta \le \alpha$ and 
$\bold q$ such that $\bold p \restriction \beta \le \bold q \in P_\beta$ 
we define
$\bold r = \bold p \otimes \bold q$ and prove that $\bold r \in P_\alpha,
\bold p \le \bold r,\bold q \le \bold r$ (in $P_\alpha$)
\item "{$(G)$}"  we prove: if for $\gamma < \beta < \alpha$ \underbar{then}
$P_\gamma < \circ P_\beta < \circ \text{ Lim}^\alpha(\bar Q)$
\item "{$(H)$}"  We define when $\bold p,\bold q \in P_\alpha$ are strongly
compatible, define their canonical common upper bound $\bold r$ and prove
$\bold r \in P_\alpha$ is a common upper bound of $\bold p,\bold q$ hence
prove strongly compatible implies compatible.
\endroster  
\endproclaim
\bigskip

\noindent
\underbar{Part (A)}:  $\bar Q = \langle P_\gamma,
{\underset\sim {}\to Q_\gamma}:\gamma < \alpha \rangle$ is from $K^{\bar a}$
if:
\medskip
\roster
\item "{$(a)$}"  each $P_\beta$ is a forcing notion
\item "{$(b)$}"  for $\beta < \alpha,{\underset\sim {}\to Q_\beta}$ is a
$P_\beta$-name of a simple forcing notion (i.e. the description
$X_{{\underset\sim {}\to Q_\beta}}$ is (see Definition \scite{1.1}) but 
for simplicity we assume that $(*)$ holds for all $\beta$ or $(*)'$ holds for
all $\beta$ where
{\roster
\itemitem{ $(*)$ }  the function $f_{{\underset\tilde {}\to Q_\beta}},
g_{{\underset\tilde {}\to Q_\beta}}$ are in $V$ (and even are objects, 
not just names) (in order to prove \scite{1.22} for one fixed countable 
field $K$ only, clearly we 
can demand that $X_{{\underset\tilde {}\to Q_\beta}} \in V$ even is an 
object not just a $P_\beta$-name but in fact the function $f,g$ in 
\scite{1.2}(3) (and $K^*,{\Cal U}^*$) were chosen so that the proof works
even for the ``for every countable $K$") and the name of 
$X_{\underset\tilde {}\to Q_\beta}$ involves $< \kappa$ many coordinates only 
(as in the names in $B(\beta)(c)$);
alternatively
\itemitem{ $(*)'$ }  more specifically 
$Q_\beta$ is $Q$ from \scite{1.4}(4) \newline
(we shall work with the $(*)'$ version).
\endroster}
\item "{$(c)$}"  for $\beta < \alpha$ we have $\bar Q \restriction \beta =
\langle P_\gamma,{\underset\sim {}\to Q_\gamma}:\gamma < \beta \rangle$
belongs to $K^{\bar a \restriction \beta}$
\item "{$(d)$}"  for $\beta < \alpha,P_\beta = \text{ Lim}^\beta(\bar Q
\restriction \beta)$.
\endroster
\bigskip

\noindent
\underbar{Part (B)}:  Let $\bar Q \in K^{\bar a}$, we let the set of elements
of $P_\alpha = \text{ Lim}^\alpha(\bar Q)$ be defined as follows.  Each
$\bold p \in P_\alpha$ will have a depth $\zeta = \zeta(\bold p) < \kappa$.
We define the set of $\bold p \in P_\alpha$ of depth $\zeta$ by induction
on $\zeta$.  Now $\bold p$ is a member of $P_\alpha$ of depth $\zeta$ if:
$\bold p = \langle p_u:u \le_e u[\bold p] \rangle$ (we may write
$p^{\bold p}_u$ or $p_u[\bold p]$ instead $p_u$ and $u^{\bold p} =$ instead
$u[\bold p]$) where:
\medskip
\roster
\item "{$(\alpha)$}"  $u[\bold p]$ is an $(\bar a,\zeta)$-creature
\item "{$(\beta)$}"  each $p_u$ (i.e. for $u \le_e u[\bold p]$) has the
form \newline
$(u,\bar s,{\bar{\bold t}}) = \langle u,\bar s^u,{\bold{\bar t}}^u \rangle
= (u,\bar s^{u,\bold p},{\bold{\bar t}}^{u,\bold p})$ such that:
{\roster
\itemitem{ (a) }  $\bar s^u = \langle s^u_\gamma:\gamma \in \text{ Dom}
(u)\rangle,s^u_\gamma$ is a function with domain \newline
$n_\gamma = \text{ Dom}(s^u_\gamma) \in \omega$ and Rang$(s^u_\gamma) 
\subseteq {\Cal H}(\aleph_0)$
\itemitem{ (b) }  ${\bold{\bar t}}^u = \langle \bold t^u_\gamma:\gamma \in
\text{ Dom}(u)\rangle$, where $\bold t^u_\gamma = 
\{{\underset\sim {}\to t^{u,\gamma}_y}:y \in Y^u_\gamma\}$; (we shall say
``$\bold t^u_\gamma$ finite" instead of ``$Y^u_\gamma$ is finite", 
similarly for $\subseteq$ etc., and demand $Y^u_\gamma \subseteq 
\{\langle \gamma,\xi,n \rangle:n < \omega,\xi < \zeta^u\}$; (clearly
$\bold t^u_\gamma = \emptyset \Leftrightarrow Y^u_\gamma = \emptyset$)
\itemitem{ (c) }  each ${\underset\sim {}\to t^{u,\gamma}_y}$ is a
$P_\gamma$-name of a member of ${}^\omega \omega$ (but $s^u_\gamma,Y^u_\gamma$
are not names!). \newline
Moreover, there is a $v = v({\underset\sim {}\to t^{u,\gamma}_y},\bold p) 
\le_e u$ (not depending on $u$) and for each $k < \omega$ there is a set 
\newline
${\Cal I} = {\Cal I}^{u,\gamma}
_{y,k} \subseteq {\Cal I}[u] =: \{ \bold r:\bold r \in P_\gamma,u[\bold r] 
\le_e v \restriction \gamma\}$, each element of ${\Cal I}$ forcing a value to
${\underset\sim {}\to t^{u,\gamma}_y}(k)$ such that every $\bold q \in
P_\gamma$ with $u \le_e u[\bold q]$ is compatible with some $\bold r \in
{\Cal I}$ and we even demand strongly compatible (see part (H)), note we use
it to $P_\gamma$; i.e. to $\bar Q \restriction \gamma$)
\itemitem{ (d) }  Moreover, if $a_\gamma$ has cardinality $\ge \kappa$, then 
$\bold t^u_\gamma$ is finite and \newline
$\bold r \in {\Cal I}^{u,\gamma}_y \Rightarrow u[\bold r] \subseteq 
a_\gamma$ (we could be more specific here).
\endroster}
\endroster
\medskip

\noindent
We further require:
\medskip
\roster
%
\item "{$(\gamma)(a)$}"  if $u_1 \le_e u_2 \le_e u[\bold p]$ \underbar{then}
\newline
$\bar s^{u_1} = \bar s^{u_2} \restriction (\text{Dom } u_1)$ (not required
from ${\bold{\bar t}}!)$ \footnote{Why the asymmetry?  Because given a
$\triangle$-system building a bound to $\aleph_0$ of them, for a specific
coordinate $\gamma$, all of the conditions contribute the same $s$ but
not necessarily the same $\bold t$.} \newline
(so $s^{u,\bold p}_\alpha$ can be written $s^{\bold p}_\alpha$)
\item "{$(b)$}"  if $u_1 \approx u_2 \le_e u[\bold p]$ \underbar{then}
\newline
${\bold{\bar t}}^{u_1} = {\bold{\bar t}}^{u_1}$
\item "{$(c)$}"  if $u_1,u_2 \le_e u[\bold p],\beta < \alpha,u_1 \restriction
\beta = u_2 \restriction \beta$ (hence $\zeta[u_1] = \zeta[u_2]$, etc.)
\underbar{then} ${\bold{\bar t}}^{u_1} \restriction \beta =
{\bold{\bar t}}^{u_2} \restriction \beta$ (for $\bar s$ this follows from the 
stronger conditions $(\gamma)(a)$).
\endroster
\medskip

$(\delta)(a)\quad$  If $u \le_e u[\bold p]$ and $\zeta^u = \xi + 1$
and $k^u_\xi = 1,v = u^{[\xi]}$ \underbar{then}:
\medskip
\roster
%
\item "{$(i)$}"  $\bar s^v = \bar s^u \restriction (\text{Dom } v)$,
(follows from $(\gamma)(a)$)
\item "{$(ii)$}"  for $\gamma \in \text{ Dom}(v)$ we have $\bold t^v_\gamma
\subseteq \bold t^u_\gamma$ and equality holds for all but finitely many
$\gamma$'s \newline
[we need the non-equal case e.g. for $t \in \bold t^u_\gamma$ whose name is
not ``based on $v$", see conditions $(\beta)(b),(c)$ above]
\item "{$(iii)$}"  for $\gamma \in \text{ Dom}(v)$ we have $\bold t^u_\gamma
\backslash \bold t^v_\gamma$ is finite
\item "{$(iv)$}" if $\gamma \in \text{ Dom}(v)$ and $a_\gamma$ has
cardinality $< \kappa$ then $\bold t^u_\gamma = \bold t^v_\gamma$
\item "{$(v)$}"  $\{\zeta:\zeta \in \text{ Dom}(u) \backslash 
\text{Dom}(u^{[\xi]}$) and
$s^u_\gamma \ne \emptyset$ or $\bold t^u_\gamma \ne \emptyset\}$ is finite
and for $\gamma \in \text{ Dom}(u) \backslash \text{Dom}(u^{[\xi]}),
Y^u_\gamma$ is finite.
\endroster
\medskip

$(b)\quad$  If $\zeta^u = \xi + 1$ and $u \le_e u[\bold p],k^u_\xi
= 2,\bar \varphi \in \Phi^{u[\bold p]}_{\zeta[u]}$ and $u = 
u[\bold p]^{[\bar \varphi]}$ and for \newline

$\qquad \qquad \quad i < \theta =: \theta^{u[\bold p]}_\xi$ 
we let $v_i = u^{(i)}$ (see Definition \scite{1.7}(10)) \underbar{then}:
\medskip
\roster
%
\item "{$(i)$}"  for $i < \theta$ and $\gamma \in \text{ Dom}(v_i) \backslash
\dsize \bigcup_{j < \theta,j \ne i} \text{ Dom}(v_j)$ we have: \newline
$\bold t^{v_i}_\gamma = \bold t^u_\gamma$ 
\item "{$(ii)$}"  for $\gamma \in \dsize \bigcap_{i < \theta} \text{ Dom}
(v_i)$ (note condition (o) of Definition \scite{1.5}) we have \newline 
$\dsize \bigcup_{i < \theta} \bold t^{v_i}_\gamma = \bold t^u_\gamma$
\item "{$(iii)$}"  for some convex equivalence relation $\bold e =
\bold e^{\bold p}_u = \bold e[\bold p,u]$ on $\theta$ with finitely many
equivalence classes we have: \newline
for $i < j < \theta$ which are $\bold e$-equivalent, $F^u_{(i,j)}$ maps
$\bold p^{[v_i]}$ to $\bold p^{[v_j]}$ \newline
(where $\bold p^{[u]},F(p^{[u]})$ are defined in Part (C) below).
\endroster
\medskip

$(c)\quad$  If $\zeta^u$ is a limit ordinal \underbar{then}
for some $\xi < \zeta^u$ we have:
\medskip
\roster
\item "{$(i)$}"  $\gamma \in \text{ Dom}(u) \backslash 
\text{ Dom}(u^{[\xi]}) \Rightarrow s^u_\gamma = \emptyset \and 
\bold t^u_\gamma = \emptyset$
\item "{$(ii)$}"  $\gamma \in \text{ Dom}(u^{[\xi]}) \Rightarrow
\bold t^u_\gamma = \bold t^{u^{[\xi]}}_\gamma$.
\endroster
\medskip

\noindent
\underbar{Part (C)}:  Let $\bold p \in P_\alpha (= \text{ Lim}^\alpha
(\bar Q))$.  For $u^* \le_e u[\bold p]$, we define \newline
$\bold p^{[u^*]} =: \langle
p_v:v \le_e u^* \rangle$ (so $u[\bold p^{[u^*]}] = u^*$).  It is easy to check
that $\bold p^{[u^*]} \in P_\alpha$.  Let \newline
$\beta \le \alpha$, we shall define
$\bold q =: \bold p \restriction \beta$ and prove that $\bold q \restriction
\beta \in P_\beta$.  If $\beta = \alpha$ let $\bold q = \bold p$, so assume
$\beta < \alpha$.  Let $\bold p = \langle p_u:u \le_e u[\bold p] \rangle,
p_u = \langle u,\bar s^u,{\bold{\bar t}}^u \rangle$; we let $u[\bold q]$ be
$u[\bold p] \restriction \beta$ (see Definition \scite{1.7}(5)) and $\bold q =
\langle q_u:u \le_e u[\bold q]\rangle$ where: \underbar{if} $u \le_e 
u[\bold q],u = v \restriction \beta$ and $v \le_e u[\bold p]$ \underbar{then} 
$q_u = \langle
u,\bar s^v \restriction \text{ Dom}(u),{\bold{\bar t}}^v \restriction
\text{ Dom}(u) \rangle$.  Now we have to prove various things.
\bigskip

\noindent
\underbar{Fact}:  $q_u$ is defined in at most one way.
\medskip
\noindent
[Why? If $v',v'' \le u[\bold p]$ and $v' \restriction \beta =
v'' \restriction \beta$ then $(v' \restriction \beta,\bar s^{v'} 
\restriction \beta,\bold{\bar t}^{v'} \restriction \beta) =$ \newline
$(v'' \restriction \beta,\bar s^{v''} \restriction \beta,\bold {\bar t}^{v''} 
\restriction \beta)$; why? first coordinate by assumption, second
coordinate by clause $(\gamma)(a)$ of Part (B) and third coordinate by
clause $(\gamma)(c)$ of Part (B)].
\bigskip

\noindent
\underbar{Fact}:  Each $q_u$ (for $u \le_e u[\bold q])$ is defined 
at least once. \newline
\medskip

\noindent
[Why? Just use \scite{1.9}(8)].
\bigskip

\noindent
\underbar{Fact}:  $\bold q \in P_\beta$.
\medskip
\noindent
[Why?  We can check all the conditions]. \newline
Let $F(\bold p) = \bold q$ if Dom$(F)$ include Dom$(u^{\bold p})$ and
$F(u^{\bold p}) = u^{\bold q}$ (see \scite{1.7}(11) and $F(p^{\bold p}_u) =
p^{\bold q}_{F(u)}$ which means: for every $\alpha \in 
\text{ Dom}(u^{\bold p})$ we have 
$s^{\bold p}_\alpha = s^{\bold q}_{F(\alpha)}$
and similarly with $\bold t$.
\bigskip

\noindent
\underbar{Part (D)}:  We define the partial orders (i.e. binary relations
which we shall prove are partial orders) $\le = \le_{nr},\le_{pr},\le_{vpr}$ 
and $\le_{apr}$ on $P_\alpha$
\medskip
\roster
\item "{$(\alpha)$}"  $\bold p^0 \le_{pr} \bold p^1$ where $\bold p^\ell =
\langle p^\ell_{\bold u}:u \le_e u[\bold p^\ell] \rangle$ (so \newline
$p^\ell_u = (u,\bar s^{u,{\bold p}^\ell},
{\bold{\bar t}}^{u,{\bold p}^\ell}$)) holds \underbar{iff} for some \newline
$u \le_e u[\bold p^1]$ and $\beta \le \alpha$ we have:
\smallskip
\noindent
{\roster
\itemitem { $(a)$ }  $u[\bold p^0] = u \restriction \beta$,
\smallskip
\noindent
\itemitem { $(b)$ }  $\dsize \bigwedge_{v \le_e u} \bar 
s^{v \restriction \beta,{\bold p}^0} =
\bar s^{v,{\bold p}^1} \restriction \beta$ and 
\smallskip
\noindent
\itemitem { $(c)$ }  $\dsize \bigwedge_{v \le_e u} 
{\bold{\bar t}}^{v \restriction \beta,{\bold p}^0} = 
{\bold{\bar t}}^{v,{\bold p}^1} \restriction \beta$
\endroster}
\smallskip
\noindent
\item "{$(\beta)$}"  $\bold p^0 \le_{vpr} \bold p^1$ if in $(\alpha),u 
\le_{de} u[\bold p^0]$ and \newline
$\gamma \in \text{ Dom}(u[\bold p^1]) \backslash \text{ Dom}(u[\bold p^0])
\Rightarrow s^{u[\bold p^1],\bold p^1}_\gamma = \emptyset \and 
\bold t^{u[\bold p^1],\bold p^1}_\gamma = \emptyset$ 
\newline
(of course by clause $(\alpha)$ we have \newline
$\beta \in u[\bold p^0] \Rightarrow
(s^{\bold p^1}_\beta,\bold t^{u[\bold p^1],\bold p^1}_\beta) =
(s^{\bold p^1}_\beta,\bold t^{u[\bold p^0],\bold p^0}_\beta)$)
\smallskip
\noindent
\item "{$(\gamma)$}"  $\bold p^0 \le_{apr} \bold p^1$ where $\bold p^\ell =
\langle p^\ell_u:u \le_e u[\bold p^\ell] \rangle$ \newline
so $p^\ell_u = (u,\bar s^{u,{\bold p}^\ell},{\bold{\bar t}}
^{u{\bold p}^\ell}))$ holds \underbar{iff}:
\smallskip
\noindent
{\roster
\itemitem{ (a) }  $u[\bold p^0] \approx u[\bold p^1]$
\smallskip
\noindent
\itemitem{ (b) }  for every $u \le_e u[\bold p^1]$ and $\gamma \in \text{ Dom}
(u)$, we have $s^{u,\bold p^0}_\gamma \subseteq s^{u,\bold p^1}_\gamma$ 
\newline
and $\bold t^{u,\bold p^0}_\gamma \subseteq \bold t^{u,\bold p^1}_\gamma$
\smallskip
\noindent
\itemitem{ (c) }  for every $\gamma \in \text{ Dom}(u[\bold p^0])$ we have
(see Definition \scite{1.1}(2)):
$$
\bold p^1 \restriction \gamma \Vdash_{P_\gamma} ``Q_\gamma \models
(s^{u[\bold p^0],\bold p^0}_\gamma,\bold t^{u[\bold p^0],\bold p^0}_\gamma)
\le (s^{u[\bold p^1],\bold p^1}_\gamma,\bold t^{u[\bold p^1],\bold p^1}
_\gamma)"
$$
\noindent
(Note: this trivially holds if $s^{u[\bold p^0],\bold p^0}_\gamma =
s^{u[\bold p^1],\bold p^1}_\gamma$ or if \newline
$s^{u[\bold p^0],p^0}_\gamma \subseteq s^{u[\bold p^1],\bold p^1}_\gamma \and
t^{u[\bold p^0],\bold p^0}_\gamma = \emptyset$
\itemitem{ (d) }  if in clause (c) the set $a_\gamma$ is of cardinality 
$< \kappa$ then we moreover have: \newline
if $\underset\sim {}\to t \in \bold t^{u[\bold p^0],\bold p^0}_\gamma$
then $(\bold p^1)^{[\underset\sim {}\to t,\bold p^0]}$ force a
value to $\underset\sim {}\to t(n)$ whenever \newline
$n \in \text{ Dom}(s^{u[\bold p^1],\bold p^1}_\gamma) \backslash 
\text{ Dom}(s^{u[\bold p^0],\bold p^0}_\gamma)$. \newline
Moreover for any such $n$, if $s^{u[\bold p^1],\bold p^1}_\gamma(n) \ne 0$
then it is computed \newline
by $g_{\underset\tilde {}\to Q_\gamma}(\in V)$ from
$s^{u[\bold p^1],\bold p^1}_\gamma \restriction n$ and the values above.
\endroster}
\item "{$(\delta)$}"  Now let $P_\alpha \models ``\bold p \le \bold r"$ iff
for some $\bold q \in P_\alpha$ we have \newline
$P_\beta \models ``\bold p \le_{pr} \bold q"$, 
and $P_\alpha \models ``\bold q \le_{apr} \bold r"$
\item "{$(\varepsilon)$}"  $\bold p^0 \approx \bold p^1$ means
$u[\bold p^0] \approx u[\bold p^1]$ \newline
and $u \le_e u[\bold p^1] \Rightarrow \bold p^0_u = \bold p^1_u$.
\endroster
\bigskip

\noindent
\underbar{Part E}:  
\roster
\item "{$(\alpha)$}"  $\bold p \le_x \bold p$ if $x \in \{nr,vpr,pr,apr\}$
\newline
[why? check]
\item "{$(\beta)$}"  if $\bold p \le_{pr} \bold q$ and $(\beta,u)$ is as in
clause $(\alpha)$ of Part (D) then $u[\bold p] \le_{rd} u[\bold q]$ hence
Dom$(u[\bold p]) \subseteq \text{ Dom}(u[\bold q])$, moreover $(\beta,u)$ is
a witnessing pair for $u[\bold p] \le_{rd} u[\bold q]$ and if $\beta$ is as
in \scite{1.10}(3), even a good witnessing pair and we say the pair witnesses
$\bold p \le_{pr} \bold q$. \newline
[why? check].
\item "{$(\gamma)$}"  if $\bold p_\ell \le_{pr} \bold q$ is witnessed by 
$(\beta,u)$ for $\ell =0,1$ then $\bold p_0 = \bold p_1$ \newline
[why? read the definitions]
\item "{$(\delta)$}"  if $\bold p_0 \le_{pr} \bold p_1 \le_{pr} \bold p_2$ 
then $\bold p_0 \le_{pr} \bold p_1$
\newline
[why?  let $(\beta_\ell,u_\ell)$ witness 
$\bold p_\ell \le_{pr} \bold p_{\ell +1}$ for $\ell = 0,1$,
so by \scite{1.10} we know some $(\beta'_0,u'_0)$ witness $u[\bold p_0]
\le_{rd} u[\bold p_2]$ and check that it witnesses $\bold p_0 \le_{pr}
\bold p_1$]
\item "{$(\varepsilon)$}"  if $\bold p_0 \le_{vpr} \bold q$ then
$\bold p_0 \le_{pr} \bold q$ \newline
[why? demanded in the definition]
\item "{$(\zeta)$}"  if $\bold p^0 \le_{vpr} \bold p^1 \le_{vpr} \bold p^2$
\underbar{then} $\bold p_0 \le_{vpr} \bold p_2$ \newline
[why? as Dom$(u[\bold p^0]) \subseteq \text{ Dom}(u[\bold p^1]) \subseteq
\text{ Dom}(u[\bold p^2])$, and if \newline
$\gamma \in \text{ Dom}(u[\bold p^2])
\backslash \text{ Dom}(u[\bold p^0])$ then for some $\ell \in \{0,1\}$ we
have \newline
$\gamma \in \text{ Dom}(u[\bold p^{\ell +1}]) \backslash
\text{ Dom}(u[\bold p^\ell])$ hence $s^{u[\bold p^2],\bold p^2}_\gamma =
s^{u[\bold p^{\ell +1}],\bold p^\ell}_\gamma = \emptyset$ \newline
(as $\bold p_\ell \le_{vpr} \bold p_{\ell +1})$ and similarly
$\bold t^{u[\bold p^2],\bold p^2}_\gamma = \bold t^{u[\bold p^{\ell +1}],
\bold p^{\ell + 1}}_\gamma = \emptyset$]
\item "{$(\eta)$}"  if $\beta < \alpha$ and $\bold p \le_* \bold q$ then in
$P_\gamma,\bold p \restriction \beta \le_x \bold q \restriction \beta$
\newline
[why?  check]
\item "{$(\theta)$}"  if $\bold p^0 \le_{apr} \bold p^1 \le_{apr} \bold p^2$
then $\bold p^0 \le_{apr} \bold p^2$ \newline
[why?  straightforward; for clause (c) note by clause $(\eta)$ that \newline
$P_\gamma \models ``\bold p^0 \restriction \gamma \le \bold p^1 \restriction
\gamma \le \bold p^2 \restriction \gamma"$]
\item "{$(\iota)$}"  if $u \le_\ell u[\bold p]$ and $\bold p \in P_\alpha$
and $\beta \le \alpha$ \underbar{then} \newline
$\bold p^{[u]} \restriction \beta \le_{pr} \bold p$ as witnessed by 
$(\beta,u)$
\item "{$(\kappa)$}"  if $\bold p^0 \le_{apr} \bold p^1 \le_{pr} \bold p^2$
\underbar{then} for some 
$\bold q$ we have $\bold p^0 \le_{pr} \bold q \le_{apr} \bold p^2$. \newline
Why?  Let $(\beta^*,u^*)$ be a good witness to $\bold p^1 \le_{pr} \bold p$.
\newline
We define $\bold q$ as follows:
{\roster
\itemitem { $(a)$ }  $u[\bold q] = u[\bold p^2]$
\smallskip
\noindent
\itemitem { $(b)$ }  $s^{u[\bold q],\bold q}$ is: 
$s^{u[\bold p^2],\bold p^2}_\gamma$  if $\gamma \in \text{ Dom}
(u(\bold p^2_2)) \backslash \text{Dom}(u[\bold p^1])$ \newline
$s^{u[\bold p^0],\bold p_0}$ if $\gamma \in \text{ Dom}(u[\bold p^1])$
\smallskip
\noindent
\itemitem { $(c)$ }  $s^{u,\bold q}_\gamma \text{ for } u \le_e u[\bold q] 
\text{ is } s^{u[\bold q],\bold q}_\gamma$
\smallskip
\noindent
\itemitem { $(d)$ }  $\bold t^{u,\bold q}_\gamma \text{ for } \gamma \in 
\text{ Dom}(u),u \le_e u[\bold q]$ is: \newline
$\bold t^{u \restriction \beta^*,\bold p^0}_\gamma \text{ \underbar{if} }
\gamma \in \text{ Dom}(u[\bold p_0]) \and u \restriction \beta^* \le_e
u[\bold p_0]$ \newline
$\bold t^{u,\bold p}_\gamma \text{ \underbar{if} otherwise}$.
\endroster}
\endroster
\medskip

\noindent
We now have to prove that $\bold q$ is as required.
\bigskip

\noindent
\underbar{$\bold q \in P_\alpha$}:  That is we have to check the definition
in Part B.

Almost all clauses of Part B are obvious, but we have to say something on
$(\delta)$(b)(iii) hence the definition of $s^{u,\bold q}_\gamma$.  So
let $u \le_e u[\bold q],\zeta^u = \xi +1,k^u_\xi = 2,u = u[\bold q]
^{[\bar \varphi]}$ where $\bar \varphi \in \Phi^{u[\bold q]}_{\zeta[u]}$ and
for $i < \theta =: \theta^{u[\bold q]}_\xi$ we let $v_i = u^{(i)}$.  So as
$u[\bold q] = u[\bold p^2]$ we can replace $\bold q$ by $\bold p_2$ and so
there is a convex equivalence relation $\bold e$ on $\theta$ such that
$i \bold e j \Rightarrow F^u_{\langle j \rangle \char 94 \bar \varphi,
\langle i \rangle \char 94 \bar \varphi}$
maps $(\bold p^2)^{[v_i]}$ to $(\bold p^2)^{[v_j]}$.  Does $\bold e$ ``work"
for the $q^{[v_i]}$?  The problem is with the $s^{v_i,\bold q}_\gamma$'s and
it works by \scite{1.9}(13).
\bigskip

\noindent
\underbar{$\bold p^0 \le_{pr} \bold q$}  Straightforward.
\bigskip

\noindent
\underbar{$\bold q \le_{apr} \bold p^2$}  Straightforward.
\bigskip

\noindent
(Compare with \scite{1.13A}(1))
\medskip
\roster
\item "{$(\lambda)$}"  $\le^{P_\alpha}$ is a partial order. \newline
[why?  follows by clause $(\kappa)$ 
as $\le^{P_\alpha}_{pr},\le^{P_\alpha}_{apr}$ are transitive]
\item "{$(\mu)$}"  if $\beta_1,\beta_2 \le \alpha,p \in P_\alpha$ then
$(\bold p \restriction \beta_1) \restriction \beta_2 = \bold p \restriction
(\text{Min}\{\beta_1,\beta_2\})$  \newline
[check]
\item "{$(\nu)$}"  $u[\bold p \restriction \beta] = u[\bold p] \restriction
\beta$ \newline
[why? check]
\item "{$(\xi)$}"  for $\beta \le \alpha,x \in \{nr,pr,apr,vpr\}$, we have
$\le^{P_\alpha}_x \restriction P_\beta = \le^{P_\beta}_x$ \newline
[why? read the definitions]
\item "{$(\circ)$}"  if $\alpha^* \le \alpha$ and $\bold p \le^{P_\alpha}
_{pr} \bold q$ \underbar{then} $(\bold p \restriction \alpha^*) \le^{P_\alpha}
_{pr} (\bold q \restriction \alpha^*)$ \newline
[why? let $\beta \le \alpha$ and $u$ be as in Part (D) clause $(\alpha)$,
witnessing $\bold p \le^{P_\alpha}_{pr} \bold q$, let $\beta^* =
\text{ Min}\{\beta,\alpha^*\}$, so it suffices to prove that $\beta^*,u
\restriction \alpha^*$ witness \newline
$(\bold p \restriction \alpha^*) \le^{P_\alpha}_{pr} 
(\bold q \restriction \alpha^*)$, and we have to check
(a),(b),(c) of $(\alpha)$ of Part (D).
\endroster
\bigskip

\noindent
\underbar{Clause (a)}:  $u[\bold p \restriction \alpha^*] = u[\bold p]
\restriction \alpha^* = (u[\bold q \restriction \beta]) \restriction
\alpha^* = u[(\bold q \restriction \beta) \restriction \alpha^*]$, \newline
$u[\bold q \restriction \text{ Min}\{\beta,\alpha^*\}] = u[(\bold q
\restriction \alpha^*) \restriction \beta^*] = u[\bold q \restriction
\alpha^*] \restriction \beta^*$.
\bigskip

\noindent
\underbar{Clause (b),(c)}:  Check.
\medskip
\roster
\item "{$(\pi)$}"  if $\alpha^* \le \alpha$ and $\bold p \le^{P_\alpha}
_{vpr} \bold q$ then $(\bold p \restriction \alpha^*) \le^{P_\alpha}
_{vpr} (\bold q \restriction \alpha^*)$ \newline
[why?  easy].
\item "{$(\rho)$}"  if $\alpha^* \le \alpha$ and $\bold p \le^{P_\alpha}
_{apr} \bold q$ then $(\bold p \restriction \alpha^*) \le^{P_\alpha}_{apr}
(\bold q \restriction \alpha^*)$ \newline
[why?  think].
\item "{$(\xi)$}"  if $\alpha^* \le \alpha$ then $\bold p \restriction
\alpha^* \le_{pr} \bold p$ \newline
[why?  check].
\endroster
\bigskip

\noindent
\underbar{Part (F)}:  If $\bold p \in P_\alpha,\beta \le \alpha,\bold p
\restriction \beta \le \bold q \in P_\beta$ \underbar{then} we shall define
$\bold r = \bold p \otimes \bold q$ and prove $\bold r \in P_\alpha,
\bold p \le \bold r,\bold q \le \bold r$.  The definition of $\bold r$ is as
follows:
\medskip
\roster
\item "{$(\alpha)$}"  $u[\bold p] \le_{de} u[\bold r] = u[\bold p] \otimes 
u[\bold q]$
\newline
(see \scite{1.10})
\item "{$(\beta)$}"  $\bar s^{u[\bold r],\bold r} 
\restriction (\text{Dom}(u[\bold p]) \backslash
\text{Dom}(u[\bold q])) \subseteq \bar s^{u[\bold p],\bold p}$
\item "{$(\gamma)$}"  $\bar s^{u[\bold r],\bold r} \restriction 
\text{ Dom}(u[\bold q]) = \bar s^{u[\bold q],\bold q}$
\item "{$(\delta)$}"  $u \le_e u[\bold r] \Rightarrow \bold t^{u,\bold r} 
\restriction \beta = \bold t^{u \restriction \beta,\bold q}$
\item "{$(\varepsilon)$}"  $u \le_e u^{\bold p} \and \gamma \in 
\text{ Dom}(u) \backslash \beta \Rightarrow
\bold t^{u,\bold r}_\gamma = \bold t^{u,\bold p}_\gamma$
\item "{$(\zeta)$}"  $u^{\bold p} \le_e u \le_e u[\bold r] \and \gamma 
\in \text{ Dom}(u) \backslash \text{Dom}(u[\bold q]) \Rightarrow 
\bold t^{u,\bold r}_\gamma = \bold t^{[u[\bold p],\bold p}_\gamma$.
\endroster
\bigskip

\noindent
Let $u^* \approx u[\bold q]$ be such that $u[\bold p \restriction \beta]
\le_{de} u^*$ (exists by \scite{1.9}(7)).  Clearly \newline
$\zeta[u[\bold p]] \le
\zeta[u^*]$; for $\xi \le \zeta[u^*]$ let $u_\xi$ be $(u^*)^{[\xi]}$.  Let
$\zeta^0 = \zeta[u[\bold p]]$ and $\zeta^1 = \zeta[u^*]$.  We now define by
induction on $\zeta \in [\zeta^0,\zeta^1]$, a condition $\bold r_\zeta$ with
$u[\bold r_\zeta] = u_\zeta$ such that $\bold r_\zeta \restriction \beta =
\langle q_u:u \le_e u_\zeta \rangle$, and $\bold r_{\zeta_0} = \bold p
\restriction \beta$.  There is no problem, (remembering $(\gamma)(b)$ of
part (B).)
\bigskip

\noindent
\underbar{Part (G)}:  For $\gamma \le \beta < \alpha,P_\gamma \le \circ
P_\beta \le \circ P_\alpha$.  Follows by parts (E),(F). 
\bigskip

\noindent
\underbar{Part (H)}:  We define: $\bold p^1,\bold p^2 \in P_\alpha$ are
strongly compatible (inside $u$) if
\medskip
\roster
\item "{$(a)$}"   for $\ell < 2$ 
we have $u^\ell \le_e u \in \text{ CR}_\alpha,
\bold p^\ell \in P_\alpha$ and $u[\bold p^\ell] = u^\ell$
\item "{$(b)$}"  if $\gamma \in \text{ Dom}(u_0) \cap \text{ Dom}(u_1)$ and
$s^{u^1,\bold p^0}_\gamma \ne s^{u^2,\bold p^1}_\gamma$ then for some
$\ell \in \{1,2\}$ we have $\bold t^{u^\ell,\bold p^\ell}_\gamma = 
\emptyset = s^{u^\ell,\bold p^\ell}_\gamma$ or at least 
$\bold t^{u^\ell,\bold p^\ell}_\gamma = \emptyset \and s^{u^\ell,\bold p^\ell}
_\gamma \subseteq s^{u^{2-\ell},\bold p^{2-\ell}}_\gamma$. 
\endroster
\medskip

\noindent
Define $\bold q$ and their canonical common upper bound by:

$$
u[\bold q] = u,\bold q = \langle q_v:v \le_e u \rangle
$$

$$
s^{v,\bold q}_\gamma = s^{u^0,\bold p^0}_\gamma \cup
s^{u^1,\bold p^1}_\gamma
$$

$$
\bold t^{v,\bold q}_\gamma = \bold t^{u^0,\bold p^0}_\gamma \cup
\bold t^{u^1,\bold p^1}_\gamma.
$$
\medskip

\noindent
\underbar{Fact}:  $\bold q$ actually is a common upper bound of 
$\bold p^0,\bold p^1$. \hfill$\square_{\scite{1.11}}$
\bigskip

\noindent
Some properties of those partial orders:
\proclaim{\stag{1.12} Claim}  1) If $\bold p \in P_\alpha$ \underbar{then} 
$[u \approx u[\bold p] \Rightarrow \bold p^{[u]} \cong \bold p]$ and
$[u \le_e u[\bold p] \Rightarrow \bold p^{[u]} \le \bold p]$ moreover
$[u \le_e u[\bold p] \Rightarrow \bold p^{[u]} \le_{pr} \bold p]$. \newline
2) If $\,\bold p^1,\bold p^2 \in P_\alpha,u[\bold p^1] \approx u[\bold p^2]$
and $[u \le_e u[\bold p^1] \Rightarrow p^1_u = p^2_u]$ \underbar{then}
${\bold{\bar p}}^1 \le_{apr} \bold p^2 \le_{apr} \bold p^1$; i.e.
${\bold{\bar p}}^1 \approx \bold p^2$. \newline
3) $(P_\alpha,\le_{vpr})$ is $\kappa$-complete, in fact any $\le_{vpr}$-
increasing sequences of length $\delta < \kappa$ has a lub. \newline
4) If $P_\alpha \models ``\bold p \le_{apr} \bold q",\beta < \alpha$ and
$\bold r$ is defined by: $u[\bold r] = u[\bold q]$, \newline
$r_u = \langle u,
\bar s^{u,\bold q} \restriction \beta \cup \bar s^{u,\bold p} \restriction
[\beta,\alpha),{\bold{\bar t}}^{u,\bold q} \restriction \beta \cup
{\bold{\bar t}}^{u,\bold p} \restriction [\beta,\alpha) \rangle$
\underbar{then} $\bold p \le_{apr} \bold r \le_{apr} \bold q$. \newline
5) If $P_\alpha \models ``\bold p \le \bold q"$ \underbar{then} for some
$\bold p^0,\bold p^1,\bold p^2 \in P_\alpha$ we have: \newline
$\bold p = \bold p^0 \le_{vpr} \bold p^1 \approx \bold p^2 \le_{apr}
\bold q$. \newline
6) If $\bold p \in P_\alpha,\zeta^{u[\bold p]}$ is a limit ordinal, $\xi <
\zeta$ large enough as in \scite{1.11}(B)$(\delta)$(c), $\bold p^{[u(\bold p)
^{[\xi]}]} \le_{apr} \bold q'$ \underbar{then} we can find $\bold q \in
P_\alpha,\bold p \le_{apr} \bold q$ such that $\bold q^{[u[\bold q]^{[\xi]}]}
= \bold q'$, and \scite{1.11}(B)$(\delta)$(c) holds for $\bold q,\xi$. 
\newline
7) If $\bold p \in P_\alpha,u \le_e u[\bold p]$, of course 
$\bold p^{[u]} \le_{pr} \bold p$ and $\bold p^{[u]} \le_{apr} \bold q \in 
P_\alpha$ \underbar{then} there is a unique
$\bold r \in P_\alpha$ such that $\bold p \le_{apr} \bold r,\bold r^{[u]} =
\bold q$ and for any $v$: \newline
$[\gamma \in \text{ Dom}(v) \backslash \text{Dom}(u[\bold q])
\and v \le_e u[\bold p] \Rightarrow s^{v,\bold p}_\gamma = s^{v,\bold r}
_\gamma \and \bold t^{v,\bold p}_\gamma = \bold t^{v,\bold r}_\gamma]$.
\newline
8)  For any $\bold p \in P_\alpha,\bar \varphi \in \Phi^{u[\bold p]}$ and 
$u^* = u[\bold p]^{[\bar \varphi^*]}$
and $\zeta < \zeta^{u[\bold p]}$ there are 
$E = E^{\bold p}_{\zeta,u^*}$ and a set $\Phi^*$ such that:
\medskip
\roster
\item "{$(a)$}"  $E$ is an equivalence relation on \newline
$\Phi^{u[\bold p]}_{\zeta,u^*} = \{ \bar \varphi \in \Phi^{u[\bold p]}_\zeta:
\bar \varphi \restriction [\zeta^{u^*},\zeta^{u[\bold p]}) = 
\bar \varphi^*\}$
\item "{$(b)$}"  $E$ has finitely many equivalence classes
\item "{$(c)$}"  if $\bar \varphi^1 E \bar \varphi^2$ then
$F^{u[\bold p]}_{\bar \varphi^1,\bar \varphi^2}$ (see \scite{1.7}(12)
is an order preserving function from $u[\bold p]^{[\bar \varphi^2]}$ onto 
$u[\bold p]^{[\bar \varphi^1]}$ 
\item "{$(d)$}"  $\Phi \subseteq \Phi^{u[\bold p]}_{\zeta,u^*}$ is a finite
set of representatives
\item "{$(e)$}"  if $\bar \varphi^1,\bar \varphi^2 \in 
\Phi^{u[\bold p]},\bar \varphi^1 E \bar \varphi^3,\bar \varphi^2 E 
\bar \varphi^4$ and $\gamma \in
\text{ Dom}(u[\bold p]^{\bar \varphi^3}) \cap \text{ Dom}(u[\bold p]
^{\bar \varphi^4})$ then \newline
$F^{u[\bold p]}_{\bar \varphi^1,\bar \varphi^3}(\gamma) = 
F^{u[\bold p]}_{\bar \varphi^2,\bar \varphi^4}(\gamma)$
\item "{$(f)$}"  if $\bar \varphi^1 E \bar \varphi^2$ \underbar{then}
$F^{u[\bold p]}_{\bar \varphi^1,\bar \varphi^2}$ maps
$\bold p^{[u^{[\bar \varphi^1]}]}$ to $\bold p^{[u^{[\bar \varphi^2]}]}$.
\item "{$(g)$}"  If $\gamma \in \text{ Dom}(u[\bold p]),\bold p 
\in P_\alpha$ and $a_\gamma$ is countable and $\underset\sim {}\to t =
t^{u[\bold p],\bold p}_y \in \bold t^{u[\bold p^0],\bold p}_\gamma$ 
\underbar{then} for
a unique $v = v[\underset\sim {}\to t,\bold p] \le_e u[\bold p]$ and
$\bar \varphi = \langle \varphi_\zeta:\xi \le \zeta < \zeta^u \rangle$ we
have $k^{u[\bold p]}_\xi = 1,\gamma \in \text{ Dom}(u^{\bar \varphi
\restriction [\xi + 1,\zeta ^{u[\bold p]}}) \backslash \text{ Dom}(u^
{\bar \varphi})$, so necessarily ${\Cal I}^{u,\gamma}_{y,k} \subseteq
{\Cal I}[u^{[\bar \varphi]}]$.  Moreover $\xi$ is the same for all
$\underset\sim {}\to t \in \bold t^{u[\bold p],\bold p}_\gamma$ (and it is
$\xi$ from \scite{1.9}(15)).
\endroster
\endproclaim
\bigskip

\demo{Proof}  E.g. \newline
8)  This is proved for each $\zeta$, by induction on $\zeta^{u[\bold p]}$.
Now the case $\zeta^{u[\bold p]} = \zeta$ is trivial and the cases
$\zeta^{u[\bold p]}$ is limit or $\zeta^{u[\bold p]} = \xi + 1 \and
k^{u[\bold p]}_\xi = 1$ follows easily by the induciton hypothesis.  So we
can assume $\zeta^{u[\bold p]} = \xi + 1 > \zeta,k^{u[\bold p]}_\xi = 2$;
let $\bold e$ be an equivalence relation on $\theta^{u[\bold p]}_\xi$ as in
clause $(\delta)$(b)(iii) of Definition \scite{1.11}; Part (B).  Let
$i_0 < \ldots < i_{n(*)-1}$ be representatives of the $\bold e$-equivalence
classes.  For $n < n(*)$ let $\bold p^n = \bold p^{[u[\bold p]^{(i_n)}]}$, so
$\zeta^{u[\bold p^n]} = \xi < \zeta^{u[\bold p]}$, so by the induction
hypothesis there is an equivalence relation $E_n$ on $\Phi^{u[\bold p^n]}
_\zeta$ as required.  Now define the $E$ we desire:
\smallskip

$\bar \varphi^1 E \bar \varphi^2$ iff $\varphi^1_\xi \bold e \varphi^2_\xi$
and $\bar \varphi^1 \restriction [\zeta,\xi) E_n \bar \varphi^2 \restriction
[\zeta,\xi)$ for each $n < n(*)$.

The checking is easy. \hfill$\square_{\scite{1.12}}$
\enddemo
\bigskip

\noindent
The following claim helps as the $Q_\beta$'s were required to be simple.
\proclaim{\stag{1.13} Claim}  If $\bar a$ is an $\alpha$-sequence and $\beta <
\alpha,\bar Q \in K^{\bar a}$ and $n < \omega$ \underbar{then} \newline
${\Cal I}^{n,0}_\beta =: \{\bold p \in P_\alpha:\beta \in \text{ Dom}(u
[\bold p])$ and $\bar s^{u[\bold p],\bold p}_\beta$ has domain including
$\{0,\dotsc,n-1\}\}$ is dense (subset of $P_\alpha$) and even \newline
${\Cal I}^{n,1}_\beta = \{\bold p \in P_\alpha:\beta \in \text{ Dom}(u
[\bold p])$ and for some $m \in [n,\omega)$ we have 
$\bar s^{u[\bold p],\bold p}_\beta(m) \ne 0\}$ is dense 
(subset of $P_\alpha$).
\endproclaim
\bigskip

\demo{Proof}  Let $\bold p \in P_\alpha$; by clauses $(\beta)(a) + 
(\delta)(a)$ of the part B of Definition \scite{1.11}, there is 
$\bold q$ such that $\bold p \le \bold q$ and $\beta
\in \text{ Dom}(u[\bold q])$.  So by a variant of \scite{1.11} Part (F)
it is enough to prove (we shall use the simplicity of 
${\underset\sim {}\to Q_\beta}$; i.e. Definition \scite{1.1}(1)):
\medskip
\roster
\item "{$(*)$}"  if $\bold p \in P_\alpha,\beta \in \text{ Dom}(u^{\bold p})$
and $\bold p \restriction \beta \le \bold q
\in P_\beta$ \underbar{then} for some $m < \omega$ for every $k < \omega$
there are $\bold r$ and $\{s^*_1,\dotsc,s^*_m\}$ such that: \newline
$\bold q \le_{apr} \bold r \in P_\beta$ and $\bold r \Vdash_{P_\beta}
``\{ {\underset\sim {}\to t} \restriction k:t \in \bold t^{u[\bold p],\bold p}
_\beta\} = \{s^*_1,\dotsc,s^*_m\}"$.
\endroster
\medskip

\noindent
If $a_\beta$ has cardinality $\ge \kappa$, this is easy (as by 
\scite{1.11} of part (B), clause $(\beta)$(c) the set 
$\bold t^{u[\bold p],\bold p}_\beta$ is finite) 
so assume $a_\beta$ has cardinality $< \kappa$. \newline
Our problem is that we do not just have to force a value to the set \newline
$\{\underset\tilde {}\to t(n):\underset\tilde {}\to t \in \bold t
^{u{[\bold p]},\bold p}_\gamma\}$ such that it is finite, but we need an 
apriory bound, one not depending
on the process.  So the natural (and first) approach induction on 
$\zeta^{\bold p}$ seems problematic.  Here our restriction in the case 
$a_\alpha$ is of cardinality $< \kappa$ (in the definition of the 
iteration, see clause $(\delta)(a)(iv)$ of Definition \scite{1.11}) help.

Let $\xi$ be as in \scite{1.9}(15) for $\gamma,\bold p$ (and also as in
\scite{1.12}(9)).

Let $\Phi^* = \{ \bar \varphi:\bar \varphi \in \Phi^{u[\bold p]}_\xi$ and
$\gamma \in \text{ Dom}(u^{\bar \varphi \restriction (\xi + 1,\zeta
^{u[\bold p]})}) \backslash \text{ Dom}(u^{\bar \varphi})\}$ and assume 
without loss of generality $(u[\bold p] \restriction \gamma) \le_{de}
u[q]$ and let $\bar \psi = \psi^{u[q]} \restriction [\zeta^{u[\bold p]},
\zeta^{u[\bold q]})$.

Let $E = E^{\bold q}_{\xi,\bar \psi}$ and $\Phi^*$ be as in \scite{1.12}(8).
Let \newline
$\Phi^\otimes = \{\bar \varphi \in \Phi^{\bold q}_\zeta:\text{ there is }
\bar \varphi' \in \Phi^* \text{ such that } \bar \varphi' \char 94
\bar \psi^* \text{ belongs to } \bar \varphi /E\}$, for \newline
$\bar \varphi = \bar \varphi' \char 94 \bar \psi \in \Phi^\otimes$ 
let $\bold t_{\bar \varphi}
= \bold t^{u[\bold p]^{\bar \varphi^1]},\bold p}$.  Now let $\Phi^1 =
\Phi^1 \cap \Phi^\otimes$ and $\bar \varphi^0,\dotsc,\bar \varphi^{n-1}$ list
\newline
$\Phi'$ and $|\bold t_{\bar \varphi}| \le n'$ for $\bar \varphi \in
\Phi^\otimes$.

Let $u^\ell = u[\bold r]^{[\bar \varphi^\ell]}$.

Now choose $m = nm'$; so to show $(*)$ let $k$ be given, 
choose $\bold q'$ such that \newline
$\bold q \le \bold q' \in P_\gamma$ and $\bold q'$ force value to
$\underset\sim {}\to t(k')$ for $k' < k,\underset\sim {}\to t \in
\dsize \bigcup_{\ell < n} \bold t^{u^\ell,\bold p}_\gamma$ moreover 
for $\ell < n,
(\bold q')^{[u^\ell]}$ force values to $\underset\sim {}\to t \restriction
k$ for $\underset\sim {}\to t \in \bold t^{u^\ell,\bold p}_\gamma$.
\medskip

For $\bar \varphi \in \Phi^\otimes$ such that $\bar \varphi E 
\bar \varphi^\ell$ let $\bold r^{\bar \varphi} = F^{\bold q}_{\bar \varphi,
\bar \varphi^\ell}(\bold r^{[u[\bold q]^{[\bar \varphi^\ell]}})$. \newline
Now we choose $\bold r$ by increasing $\bold q^{[u]}$ if $u \in \{u[\bold q]
^{\bar \varphi}:\bar \varphi \in \Phi^\otimes\}$ and essentially only there.
Explicitly:
\medskip
\roster
%
\item "{$(i)$}"  $u[\bold r] = u[\bold q]$
\item "{$(ii)$}"  if $\beta \in \text{ Dom}(u[\bold q]^{[\bar \varphi]})$
for some $\bar \varphi \in \Phi^\otimes$ then
$$
s^{u[\bold q],\bold r}_\beta = s^{u[\bold r^{\bar \varphi}],\bold r^
{\bar \varphi}}_\gamma
$$
(well defined by \scite{1.12}(8) clause $(c)$
\item "{$(iii)$}"  if $\beta \in \text{ Dom}(u[\bold q])$ and clause (ii) 
does not apply then $s^{u[\bold r],\bold r}_\beta$ \newline
is chosen as $s^{u[\bold q],\bold q}_\gamma$
\item "{$(iv)$}"  if $\beta \in \text{ Dom}(u),u \le_e u[\bold q]
^{[\bar \varphi]},\bar \varphi \in \Phi^\otimes$ \underbar{then}
$\bold t^{u,\bold r} = \bold t^{u,\bold r^{\bar \varphi}}$ \newline
(again by \scite{1.12}(8) it is well defined)
\item "{$(v)$}"  if $\beta \in \text{ Dom}(u),u \le u[\bold q]$ and clause
(iv) does not apply then

$$
\bold t^{u,\bold r}_\beta = \bold t^{u,\bold q}_\beta \cup
\dsize \bigcup\{\bold t^{v,\bold q}_\beta:v <_e u\}
$$
\noindent
(i.e. the definition is by induction on $\zeta^u$).
\endroster
\medskip

\noindent
We leave the rest to the reader. \hfill$\square_{\scite{1.13}}$
\enddemo
\bigskip

\demo{\stag{1.13A} Fact}  1)  If $u_1 \le_e u_2 \in \text{ CR}_\alpha,
\bold p_1 \in P_\alpha$ and $u[\bold p_1] = u_1$ \underbar{then} there is
$\bold p_2,\bold p_1 \le_{vpr} \bold p_2 \in P_\alpha$ such that
$u[\bold p_2] = u_2$ in fact if $\bold p_1 \le \bold p'_2$ and 
$u[\bold p'_2] = u_2$ then $\bold p_2 \le_{apr} \bold p'_2$. \newline
2)  Suppose $u^0 \le_{de} u^1 \le_{de} u^2,u[\bold p^0] = u^0,
u[\bold p^2] = u^2$ and $\bold p^0,\bold p^1$ are strongly compatible 
\underbar{then} $(\bold p^2)^{[u^1]},\bold p^1$ are strongly compatible. 
\newline
In fact, if $\bold p_0,\bold p_2$ are strongly compatible in $u$ and
$p'_\ell = \bold p^{[u_\ell]}_\ell$ then $\bold p^1_0,\bold p'_1$ are strongly
compatible.
\enddemo
\bigskip

\demo{Proof}  1) Check the definitions. \newline
2) Check the definitions. \newline
\enddemo
\bigskip

\proclaim{\stag{1.14} Claim}  1) For any ordinal 
$\alpha$ and $\alpha$-sequence $\bar a$, for any $\bar Q \in K^{\bar a}$ 
and \newline
$u^* \in \text{ CR}_{\bar a}$,
the partial order $(\{ \bold p \in P_\alpha:
u[\bold p] \le_e u^*\},\le)$ satisfies the c.c.c. and even the Knaster 
condition. \newline
2)  If $\kappa = \aleph_1$ and $u \in \text{ CR}_\alpha$ and
$\bar Q \in K^\alpha$ \underbar{then} the set
$\{ \bar s^{u,\bold p}:\bold p \in P_\alpha,u = u[\bold p]\}$ is countable.   
\newline
3) If $\bold p_\sigma \in P_\alpha,u(\bold p_\sigma) = u^* (\in \text{ CR}
_{\bar a}),\sigma < \omega_1$ \underbar{then} for some uncountable 
$A \subseteq \omega_1$, for any $\sigma_1,\sigma_2 \in A$ the 
conditions $\bold p_{\sigma_1},\bold p_{\sigma_2}$ are strongly compatible
(in $u^*$) and so have a (canonical) common upper bound. \newline
4) For any $\bar Q \in K^\alpha$ and $\bar u = \langle u_\zeta:\zeta <
\omega_1 \rangle$ such that $u_\zeta \in \text{ CR}_{\alpha,\zeta}$ and
$\bar u$ is \newline
$\le_{de}$-increasing continuous \underbar{we have}:
$(\{\bold p \in P_\alpha:u[\bold p] \le u_\zeta$ for some $\zeta < \omega_1\},
\le)$ satisfies the c.c.c. 
\endproclaim
\bigskip

\demo{Proof}  1)  So it follows from part (3) by \scite{1.13A}(1). \newline
2)  By induction on $\zeta^*$ - (essentially included in the proof of 
part (3).  If we are satisfied with the case
$\kappa = \aleph_1$, this simplifies the proof of parts (3),(4).) \newline
3)  We prove by induction on $\zeta^{u^*}$. 
\enddemo
\bigskip

\noindent
\underbar{First Case}:  $\zeta^{u^*} = 0$. \newline
Trivial.
\bigskip

\noindent
\underbar{Second Case}:  $\zeta^{u^*} = \xi + 1,k^{u^*}_\xi = 1$.  Let
$v^* = (u^*)^{[\xi]},\bold p'_\sigma = \bold p^{[v^*]}_\sigma$.  By the
induction hypothesis without loss of generality the conclusion of (3)
holds for \newline
$\langle \bold p'_\sigma:\sigma < \omega_1 \rangle,A = \omega_1$.
By Definition \scite{1.11} part(B),clause$(\beta)(e)$ if $\kappa =
\aleph_1$ the set
$\{ \bar s^{u^*,\bold p_\sigma} \restriction (\text{Dom}(u^*) \backslash
\text{Dom}(v^*):\sigma < \omega_1\}$ is countable so without loss of
generality \newline
$\bar s^{u^*,\bold p_\sigma} \restriction (\text{Dom}(u^*) \backslash
\text{Dom}(v^*))$ for $\sigma \in A$ is constant.  
But even not assuming $\kappa = \aleph_1$,
letting $a_\sigma = \{\gamma:\gamma \in \text{ Dom}(u^*) \text{ but }
\gamma \notin \text{ Dom}(v^*) \text{ and } s^{u^*,\bold p_\sigma}_\gamma \ne
\emptyset \and \bold t^{u^*,\bold p_\sigma}_\gamma \ne \emptyset\}$ 
is finite hence without loss of generality $\langle a_\sigma:
\sigma < \omega_1 \rangle$ form a $\triangle$-system with heart say $a^*$ and
$\bar s^{u^*,p_\sigma} \restriction a^* = \bar s^*$; i.e. is constant.
Now clearly $\langle \bold p_\sigma:\sigma < \omega_1 \rangle$ are 
pairwise strongly compatible as required. 
\bigskip

\noindent
\underbar{Third Case}:   $\zeta^{u^*} = \xi + 1,k^{u^*}_\xi = 2$.  For each
$\sigma < \omega_1$ we have an equivalence relation $\bold e_\sigma$ on
$\theta^{u^*}_\xi$, as in Definition \scite{1.11} part (B), clause 
$(\gamma)(f)(iii)$, so it is convex with finitely many equivalence classes,
so the number of possible $\bold e$'s is $\le | \dsize \sum_n \,
\theta^{u^*}_\xi|^n$; as
$\theta^{u^*}_\xi < \omega_1$ without loss of generality $\bold e_\sigma 
= \bold e$.  Let $a \subseteq \theta^{u^*}_\xi$ be a set of representatives.
Let $v_i = (u^*)^{(i)}$ for $i <
\theta^{u^*}_\xi$.  By successive application of the induction hypothesis to
$\langle \bold p^{[v_i]}_\sigma:\sigma < \omega_1 \rangle$ without loss of 
generality for $i \in a$ we have $\langle \bold p^{(v_i)}_\sigma:\sigma
< \omega_1 \rangle$ satisfies the conclusion for $A = \omega_1$.
\bigskip

By use of the $F^{u^*}_{(i,j)}$ this holds for any $i < \theta^{u^*}_\xi$.
Now it suffices to prove that for any $\sigma_1,\sigma_2 < \omega_1$, the
conditions $\bold p_{\sigma_1},\bold p_{\sigma_2}$ are strongly compatible;
so let $\gamma \in u[\bold p_{\sigma_1}] = u[\bold p_{\sigma_2}] = u^*$.
If $\gamma \notin \dsize \bigcap_{i < \theta_\xi} \text{ Dom}(v_i)$ then
(see Definition \scite{1.5}) for exactly on $i = i(\gamma)$ we have
$\gamma \in v_i$, so necessarily $(s^{u^*,\bold p_{\sigma_\ell}}_\gamma,
\bold t^{u^*,\bold p_{\sigma_\ell}}_\gamma) = (s^{v_i,
\bold p_{\sigma_\ell}}_\gamma,\bold t^{v_i,\bold p_{\sigma_\ell}}_\gamma)$ 
for $\ell = 1,2$ and by the
application of the induction hypothesis above we get the desired conclusion.
So assume $\gamma \in \dsize \bigcap_{i < \theta_\xi} \text{ Dom}(v_i)$;
now if $a_\gamma$ has cardinality $< \kappa$, then by Definition \scite{1.11},
Part (B) clause $(\beta)(d)$ and $(*)$ of Definition \scite{1.8}(2) for all 
$i < \theta_\xi$, we get the same pair
$(s^{v_i,\bold p_{\sigma_\ell}}_\gamma,\bold t^{v_i,\bold p_{\sigma_\ell}}
_\gamma)$, so we are done, so assume also $a_\gamma$ has cardinality
$< \kappa$.

Also if $s^{v_i,\bold p_{\sigma_1}}_\gamma = s_\gamma^{v_i,\bold p
_{\sigma_2}}$ we are done, so by symmetry without loss of generality
$s^{v_i,\bold p_{\sigma_1}}_\gamma \subset s^{v_i,\bold p_{\sigma_2}}_\sigma$
but the $s$'s do not depend on $i$ and by what we do for each $i$ (and
\scite{1.13A}) we have $\bold t^{v_i,\bold p_{\sigma_1}}_\gamma =
\emptyset$, hence $\bold t^{u[\bold p_{\sigma_1}],\bold p_{\sigma_1}}_\gamma =
\dsize \bigcup_{i < \theta_\xi} \bold t^{u[\bold p_{\sigma_1}],
\bold p_{\sigma_1}}_\gamma = \emptyset$, and we are done. 
\bigskip

\noindent
\underbar{Fourth Case}:  $\zeta^{u^*}$ limit and we know (see Definition
\scite{1.11} part(B), clause $(\beta)(d)$) that for each 
$\sigma < \omega_1$, there
is $\zeta_\sigma < \zeta^{u^*}$ such that $[\gamma \in \text{ Dom}(u^*)
\backslash \text{Dom}(u^{[\xi]}) \Rightarrow s^u_\gamma = \emptyset =
\bold t^u_\gamma]$ (i.e. we only duplicate).  If $\kappa = \aleph_1$, 
clearly without loss of
generality $\zeta_\sigma = \zeta^*$ for $\sigma < \omega_1$; now we use
the induction hypothesis on $\langle p^{[\zeta^*]}_\sigma:\sigma <
\omega_1 \rangle$ using \scite{1.12}(6).  So assume $\kappa > \aleph_1$ and
let $\langle \varepsilon_\sigma:\sigma < \text{ cf}(\zeta^{u^*}) \rangle$
be increasing continuous with limit $\zeta$.  If cf$(\zeta) = \aleph_0$, then
for some $\sigma(*) < \omega_1,A =: \{\sigma < \omega_1:\zeta_\sigma <
\varepsilon_{\sigma(*)}\}$ is uncountable and continue as above with
$\zeta^* = \varepsilon_{\sigma(*)}$.  If cf$(\zeta) > \aleph_1$, for some
$\zeta^* < \zeta$ we have $(\forall \sigma < \omega_1)[\zeta_\sigma < 
\zeta^*]$ and continue as above.  We are left with the case cf$(\zeta) =
\aleph_1$, which has the same proof as part (4) (using the induction
hypothesis.

Now for each limit $\sigma < \omega_1$, there is $\xi_\sigma < \varepsilon
_\sigma$ such that $\alpha \in \text{ Dom}((u^*)^{[\varepsilon_\sigma]})
\and \alpha \notin \text{ Dom}((u^*)^{[\xi_\sigma]}) \Rightarrow
s^{u^*,\bold p_\sigma}_\alpha = \emptyset$.  By Fodor's lemma for some
stationary $s \subseteq \omega_1$, we have $\sigma \in S \Rightarrow
\xi_\sigma = \xi^*$, and without loss of generality $\sigma_1 < \sigma_2$
in $s \Rightarrow \zeta_{\sigma_1} < \varepsilon_{\sigma_2}$.  We now apply
the induction hypothesis to $\langle p'_\sigma - p^{[(u^*)^{[\xi^*]}]}:
\sigma \in s \rangle$. \newline
4)?  Like the proof of the last possibility in the proof of part (3)
(cf$(\zeta^u) = \aleph_1)$. \newline
5) If $u \in CR_\alpha$ and $\zeta \le \zeta^u$ then there is an
equivalence relation $E$ on $\Phi^u_\zeta$ with finitely many equivalence
classes.
\bigskip

For each $\sigma < \omega_1$ choose $\zeta_\sigma < \omega)1$ such that
$u^{\bold p_\sigma} \le_e u_{\zeta_\sigma}$ so without loss of generality
$\zeta_\sigma > \sigma$ choose $\bold p'_\sigma \in P_\alpha$ such that
$u[\bold p'_\sigma] \approx u_{\zeta_\sigma}$ and $\bold p_\sigma \le_{vpr}
\bold p'_\sigma$ and let $\bold p''_\sigma \approx \bold p'_\sigma,
u[\bold p''_\sigma] = u_{\zeta_\sigma}$, and for limit $\sigma < \omega_1$
choose $\xi_\sigma < \zeta[u_\sigma]$ such that $\gamma \in \text{ Dom}
\left( (u_{\zeta_\sigma})^{[\sigma]} \right) \and \gamma \notin \text{ Dom}
\left( (u_{\zeta_\sigma})^{[\xi_\sigma]} \right) \Rightarrow 
s^{u_{\zeta_\sigma},{\bold p''}_\sigma}_\gamma = \emptyset$ without loss of
generality $\xi_\sigma = \zeta[u_{\zeta'_\sigma}],\zeta'_\sigma < \sigma$ and
choose a stationary set of limit ordinals $< \omega_1$ such that $\sigma \in
S \Rightarrow \zeta''_\sigma = \zeta^*$ and $\sigma_1 < \sigma_2$ in
$S \Rightarrow \zeta^*_{\sigma_1} < \sigma_2$.  Let $v =$?

Apply part (3) on $\langle (p''_\sigma)^{[u_{\zeta^*}]}:\sigma \in S \rangle$
get an uncountable $S_2 \subseteq S$; now \newline
$\langle p_\sigma:\sigma \in
S_1 \rangle$ is as required.  \hfill$\square_{\scite{1.14}}$ 
\bigskip

\proclaim{\stag{1.15} Claim}  Let $\bold p \in P_\alpha$ and let
${\underset\sim {}\to \tau}$ be a $P_\alpha$-name of an ordinal.  
\underbar{Then} there is $\bold q,\bold p \le_{vpr} \bold q \in P_\alpha$ 
such that:
\medskip
\roster
\item "{$(*)$}"  if $\bold q \le \bold r \in P_\alpha,\bold r \Vdash
``{\underset\sim {}\to \tau} = \sigma"$ \underbar{then} for some
$\bold r'$, letting $\bar \varphi$ be such that $u[\bold q] =
u[\bold r']^{[\bar \varphi]},\bold q' = (\bold r')^{[\bar \varphi]}$
we have $\bold q' \Vdash ``{\underset\sim {}\to \tau} = \sigma"$.
\endroster
\medskip

\noindent
Moreover
\roster
\item "{$(**)$}"  there is a countable family \newline
${\Cal I} \subseteq \{\bold r':\text{for some } u,u^{\bold p} \le_{de}
\bold r' \le_{de} \bold q$, and $\bold q^{[u]} \le_{apr} \bold r'$ and
$\bold r'$ force a value to $\tau\}$ such that $q \le \bold q' \and
(\bold q'$ force a value to $\underset\sim {}\to \tau) \Rightarrow \bold q'$
is strongly compatible with some $\bold r' \in {\Cal I}$.
\endroster
\endproclaim
\bigskip

\demo{Proof}  Assume that there is no $\bold q$ satisfying $(**)$.
We choose by induction on $\varepsilon < \omega_1,\bold q_\varepsilon,
\bold r_\varepsilon$ such that:
\medskip
\roster
\item "{$(a)$}"  $\langle \bold q_\varepsilon:\varepsilon < \omega_1 \rangle$
is $\le_{vpr}$-increasing (in $P_\alpha$) \newline
(so $\langle u[\bold q_\varepsilon]:\varepsilon < \omega_1 \rangle$ is 
$\le_{de}$-increasing continuous)
\item "{$(b)$}"  for each successor 
$\varepsilon$ there is $\bold r_\varepsilon,\bold q_\varepsilon \le_{apr}
\bold r_\varepsilon,\bold r_\varepsilon$ forces a value to 
$\underset\sim {}\to \tau$, and
$\bold r_\varepsilon$ is not strongly compatible with any 
$\bold r_\zeta,\zeta < \varepsilon$ successor.
\endroster
\medskip

\noindent
For $\varepsilon = 0$ let $\bold q_0 = \bold p$, for limit $\varepsilon$ use
\scite{1.12}(3), for $\varepsilon$ successor use the assumption that ``there 
is no $\bold q$ as required"; i.e. suppose we cannot find such
$\bold r_\varepsilon$.

Let ${\Cal I}_\varepsilon = \{\bold r:\bold q_{\varepsilon-1} \le_{apr}
\bold r$ and $\bold r$ force a value to $\underset\sim {}\to t\}, 
\bold r_{\zeta +1}:\zeta + 1 < \varepsilon$ (so if
$\bold q_{\varepsilon - 1} \le \bold r$ and force a value to 
$\underset\sim {}\to t$ hence for some
successor ordinal $\zeta < \varepsilon$ we have $\bold r,\bold r_\zeta$ are
strongly compatible.  

So we are assuming there is no $\bold q$ as required in $(*)$ we can carry
the induction.

Lastly, $\langle u[\bold q_\varepsilon]:\varepsilon < \omega
_1 \rangle$ and $\langle \bold r_\varepsilon:\varepsilon < \omega_1 \rangle$
contradict \scite{1.14}(4). \newline
Why $(**) \Rightarrow (*)$?  So let $\bold q_{\varepsilon -1},I_\varepsilon
= \{\bold r_{\zeta + 1}:\zeta +1 < \varepsilon\}$, be as gotten above.  Now
if $\bold q_{\varepsilon-1} \le \bold r \in P_\alpha$, then for some
$\bold r',\bold r \le \bold r' \in \bold p_\alpha,\bold r'$ force a value
to $\underset\sim {}\to \tau$.  So by \scite{1.13A}(3), 
$\bold r^{[u_{\varepsilon -1}]},
\bold r_\zeta$ are strongly compatible and so their canonical common upper
bound is in ${\Cal I}_\varepsilon$ and easily is $\le$ the canonical common
upper bound of $\bold r,\bold r_\zeta$, hence ${\Cal I}_\varepsilon$
satisfies: for a dense set of $\bold r \in P_\alpha,u[\bold q_{\varepsilon
-1}] \le_e \bold r \Rightarrow \bold r^{[u[\bold q_{\varepsilon -1}]]} \in
{\Cal I}_\varepsilon$; as required.  \hfill$\square_{\scite{1.15}}$
\enddemo
\bigskip

\noindent
\underbar{\stag{1.16} Conclusion}  Forcing with $P_\alpha$ (for $\bar Q \in
K^\alpha$) does not collapse $\aleph_1$, in fact the regularity of any
cardinal $\le \kappa$ is preserved; and the forcing notion is proper.
\bigskip

\demo{Proof}  By 1.14, 1.15 (using 1.5).
\enddemo
\bigskip

\proclaim{\stag{1.17} Claim}  Let $\bar Q \in K^{\bar a}$.  \newline
1)  For any $\beta < \alpha$ the set $\{ \bold p \in P_\alpha:\beta \in
u[\bold p]\}$ is dense (and open). \newline
2)  For any $\beta < \alpha$ satisfying 
``$a_\beta$ has cardinality $< \kappa$" and 
$P_\beta$-name ${\underset\sim {}\to \tau}$ of a member of 
$\,{}^\omega \omega$ defined by
elements with support $\subseteq a_\beta$ the set
$\{\bold p_\alpha \in P_\alpha:{\underset\sim {}\to \tau} \in
\bold t^{u[\bold p],\bold p}_\beta\}$ is a dense subset of $P_\alpha$.
\endproclaim
\bigskip

\demo{Proof}  1) We can use the case $k_\xi = 1$, i.e. for any $\bold p \in
P_\alpha$, such that $\beta \notin u^* u[\bold p]$ define $v^\otimes$ by:
$v^\otimes \in \text{ CR}_{\alpha,\zeta,u[\bold p]},(v^\otimes) = u^*$, 
Dom$(v^\otimes) = \text{ Dom}(u^*) \cup \{ \beta\}$ and 
$\bold q = \langle q_v:v \le_e v^\otimes \rangle,q_v = p_v$ if $v \le_e 
u[\bold p]$ and $s^{v^*,\bold q}_\gamma$ is $s^{u^*,\bold p}_\gamma$
if $\gamma \in u^*$ and $\emptyset$ otherwise, $\bold t^{v,\bold q}_\gamma =
\bold t^{v,\bold p}_\gamma$ if $v \le_e u^*$ and $\emptyset$ otherwise.
Easily $\bold p \le \bold q,\beta \in \text{ Dom}(u[\bold q])$.
\hfill$\square_{\scite{1.17}}$
\enddemo
\bigskip

\proclaim{\stag{1.18} Conclusion}  1) If $\bar Q \in K_{\bar a}$ and 
$a_\beta$ has cardinality $\ge \kappa$, \underbar{then} \newline
$\Vdash_{P_\alpha} ``{\underset\sim {}\to s_\alpha} = \cup
\{s^{u[\bold p],\bold p}_\alpha:\bold p \in
{\underset\sim {}\to G_{P_{\alpha + 1}}},\alpha \in \text{ Dom}(\bold p)\}$
is in ${}^\omega \omega$, every initial segment is in 
$S_{{\underset\sim {}\to Q_\alpha}}$ (see Definition \scite{1.1})), 
for infinitely
many $n,{\underset\sim {}\to s_\alpha}(n) \ne 0$ and for every finite
$\bold t \subseteq ({}^\omega \omega)^{V^{P_\alpha}}$ defined by
$P_\alpha$-names described by conditions of $P_\alpha$ with support
$\subseteq a_\alpha$, for some $n$, for every
$m > n$ we have $({\underset\sim {}\to s_\alpha} \restriction n,\bold t) \le
({\underset\sim {}\to s_\alpha} \restriction m,\bold t)"$.
\endproclaim
\bigskip

\proclaim{\stag{1.19} Claim}  Assume $\kappa$ is regular uncountable and 
$\bar a$ is a $(\mu,\alpha)$-sequence.
Assume further $\bar Q \in K^a$ and $P_\alpha = \text{ Lim}^{\alpha,\bar a}
\bar Q$. \newline
1) $P_\alpha$ satisfies the $(2^{< \kappa})^+$-c.c.c. (if $2 \in \Theta$).
\newline
2)  Assume $\mu = \text{ cf}(\mu) > |\alpha|^{< \kappa}$ for $\alpha <
\kappa$ (so $\mu > \kappa$) and $\gamma < \alpha \Rightarrow |a_\gamma| <
\mu$. 
If $S \subseteq \{ \delta < \mu:\text{cf}(\delta) = \aleph_0\}$ is
stationary and $\clubsuit_S$ holds in $V$ say exemplifies by
$\bar c = \langle c_\delta:\delta \in S \rangle$ (and $\omega \in \Theta$), 
\underbar{then} forcing with $P_\alpha$ preserves it, even for the same
$\bar c$. \newline
3) Assume cf$(\alpha) > \mu$ and $\alpha = \sup\{\beta:a_\beta$ has 
cardinality $< \kappa\}$
and $\bold p \in P_\alpha$ and $\bold p \Vdash 
``{\underset\sim {}\to f_i} \in {}^\omega \omega$ for $i < \mu"$.
\underbar{Then} for some $\bold q$ and $\beta$ we have: $\bold p \le
\bold q \in P_\alpha$ and $\bold q \Vdash_{P_\alpha}$ 
``for unboundedly many $i < \mu$ we have
${\underset\sim {}\to f_i} \in \bold t^{u[\bold r],r}_\beta$ for some
$\bold r \in {\underset\sim {}\to G_{P_\alpha}}$; hence for any finite set
$w$ of such $i$'s, some $n$ for every $m > n$ we have
$({\underset\sim {}\to s_\beta} \restriction n,\{ {\underset\sim {}\to f_i}:
i \in w\}) \le ({\underset\sim {}\to s_\beta} \restriction n,\{
{\underset\sim {}\to f_i}:i \in w\})$".
\endproclaim
\bigskip

\remark{\stag{1.20} Remark}  1)  We are most interested in the case 
otp$(c) = \omega$, but any countable ordinal is O.K.
\endremark
\bigskip

\demo{Proof}  1) Let $\bold p_i \in P_\alpha$ for $i < (2^{< \kappa})^+$.  
We can define
the set of isomorphism types: it has cardinality $\le 2^{< \kappa}$ as all
${\underset\sim {}\to t^{u,\bold p_i}_{\gamma,y}}$ is a name involving only
$< \kappa$ conditions.  So without loss
of generality this isomorphism type is fixed.  Let $A_i = \cup\{a_\beta:
\text{for some } j < i,\beta \in u^{\bold p_j}\}$, so $|A_i| < \kappa$, and
$A_i$ is increasing continuous in $i$.  Without loss of generality, for some
stationary $A \subseteq \{\delta < (2^{< \kappa})^+:\text{cf}(\delta) = 
\kappa\}$ we have $\langle u[\bold p_i]:i \in A \rangle$ is a 
$\triangle$-system; i.e. $\langle \text{Dom}(u^{\bold p_i}):i \in A \rangle$ 
is a $\triangle$-system, otp(Dom$(u^{\bold p_i})$ for $i \in A$ is constant
and letting $F_{i,j}$ be the (one to one) order preserving function from
Dom$(u^{\bold p_j})$ onto Dom$(u^{\bold p_i})$, so for $i,j \in A,F_{i,j}$
maps $u^{\bold p_i}$ to $u^{\bold p_j}$ and $\bold p_i$ to $\bold p_j$.

In particular $\bar s^{u[\bold p_i],\bold p_i},\bar s^{u[\bold p_j],
\bold p_j}$ agree on the intersection of their domains.
So for $i,j \in A$ (as $2 \in \Theta$) there is $u \in \text{ CR}_{\bar a},
u_i \le_e u \and u_j \le_e u$, clearly $\bold p_i,\bold p_j$ are strongly
compatible inside $u$ (see \scite{1.13A}(2)) hence are compatible, as
required.  \newline
2)  Let $\bold p \in P_\alpha,\bold p \Vdash_{P_\alpha} 
``{\underset\sim {}\to \tau} \subseteq \mu = \sup 
{\underset\sim {}\to \tau}"$.  For each $i < \mu$ let $\bold p_i,j_i$
be such that: $i < j_i < \mu,\bold p \le \bold p_i$ and $\bold p_i \Vdash
``j_i \in {\underset\sim {}\to \tau}"$ and without loss of generality
$i_1 < i_2 \Rightarrow j_{i_1} < j_{i_2}$.  We can find an unbounded
$A \subseteq \mu$ such that the sequence $\langle p_i:i \in A \rangle$ is as
in the proof of part (1).  In the proof of part (1) without loss of 
generality; by possibly thinning $A$
(preserving its cardinality even its being stationary we can have: if
$\gamma \in \dsize \bigcap_{i \in A} \text{ Dom}(u^{\bold p_i})$ and
$|a_\gamma| < \mu$ then $\langle a_\gamma \cap \text{ Dom}u[\bold p_i]:i \in
A \rangle$ is constant.  Now any $\omega$ member of $\{\bold p_i:i \in
A\}$ has a common upper bound. 
Now assume 
$\langle a_\delta:\delta \in S \rangle (a_\delta \subseteq \delta = 
\sup(a_\delta$) and otp$(a_\delta) \in \Theta)$ exemplifies $\clubsuit_S$ 
in $V$. \newline
For some $\delta,a_\delta \subseteq \{j_i:i \in A\}$, and we can finish.
\newline
3) As cf$(\alpha) > \mu > 2^{< \kappa}$ and as $P_\alpha$ satisfies the 
$(2^{< \kappa})^+$-c.c.c. and
as $\langle P_\beta:\beta < \alpha \rangle$ is $\le \circ$-increasing 
continuous in limits of cofinality $\ge \kappa$ with limit $P_\alpha$, we 
know that for some $\beta^* <
\alpha, \langle {\underset\sim {}\to f_i}:i < \kappa \rangle$ is a
clearly $P_{\beta^*}$-name, $\bold p \in P_{\beta^*}$.  Choose 
$\gamma \in (\beta^*,
\alpha)$ such that $a_\gamma = \emptyset$, so easily for each $i$ there is
$\bold q_i$ satisfying $\bold p \le \bold q_i \in P_\alpha,\gamma \in
u[\bold q_i]$ and $\bold t^{u,\bold q_i}_\gamma = 
\{{\underset\sim {}\to f_i}\}$.  As $P_\alpha$ satisfies the 
$(2^{< \kappa})^+$-c.c.c. hene the $\mu$-c.c.c. clearly for every 
$i < \mu$ large enough we have:

$$
\bold q_i \Vdash_{P_\alpha} ``\{j < \mu:\bold q_j \in
{\underset\sim {}\to G_{P_\alpha}}\} \text{ is unbounded in } \mu".
$$
\medskip

\noindent
This clearly suffices.  \hfill$\square_{\scite{1.19}}$
\enddemo
\bigskip

\demo{\stag{1.21} Conclusion}  $(V \models CH)$.  
For some $(\aleph_1,\omega_3)$-sequence $\bar a$ and $\bar Q \in K^{\bar a},
P_{\omega_3}$ is a forcing of cardinality $\aleph_3$, satisfies the
$\aleph_2$-c.c.c., is proper and
\medskip
\roster
\item "{$(a)$}"  forcing with $P_{\omega_3}$ preserves $\clubsuit_S$ for
$S \subseteq S^2_0 =: \{\delta < \aleph_2:\text{cf}(\delta) = \aleph_0\}$
\item "{$(b)$}"  for every simple $Q$ with $g,f \in V$ and $B \subseteq
{}^\omega \omega,|B| < \aleph_2$ there is $s \in {}^\omega \omega$ such that
for every finite $\bold t \subseteq B$ for every large enough $n,Q \models
(s \restriction n,\bold t) \le (s \restriction (n+1),\bold t)$ in particular
\item "{$(b)'$}"  in $V^{P_{\aleph_3}}$, if $F$ is a countable fixed
${\Cal U}$ a vector space over $F$ of dimension $\aleph_0$, and $f_i \in
{\Cal U}_F$ a homomorphism for $i < \aleph_1$ \underbar{then} there is a
independent sequence $\langle x_i:i < \omega \rangle$ of elements of
${\Cal U}$ such that \newline
$(\forall i < \omega_1)$ [for every $n$ large enough, $f_i(x_n) = 0$]
\item "{$(c)$}"  if $F,{\Cal U}$ is as above, $f_i \in {}^{\Cal U}F$ a
homomorphism for $i < \omega_2$ then for some unbounded $B \subseteq
{}^\omega 2$ and an independent sequence $\langle x_n:n < \omega \rangle$ of
members of ${\Cal U}$ \newline
$(\forall i \in B)$ [for every $n$ large enough, $f_i(x_n) = 0$].
\endroster
\enddemo
\bigskip

\demo{Proof}  We know cf$([\aleph_3]^{\aleph_1},\subseteq) = 
\aleph_3$ hence we can find
$\bar a = \langle a_\alpha:\alpha < \aleph_3 \rangle$ such that \newline
$a_\alpha \in [\alpha]^{\le \aleph_1}$ and $(\forall a)(a \subseteq 
\aleph_3 \and
|a| \le \aleph_1 \rightarrow (\exists^{\aleph_3} \beta < \aleph_3)(a
\subseteq a_\beta))$ and $(\forall a)(a \subseteq \aleph_\alpha \and |a| = 
\aleph_0 \rightarrow (\exists^{\aleph_3} \beta < \aleph_3)(a = a_\beta))$.
Now choose $\bar Q \in K^{\bar a}$ such that for every $\gamma
< \aleph_3$ and $P_\gamma$-name of a simple forcing notion $Q$ with
$g,f \in V$ the set $\{\beta:Q_\beta = Q[V^{P_\beta}]\}$ is
unbounded in $\aleph_3$.  Alternatively deal only with $Q$ of
\scite{1.3}(4) and derive the results for the others.  
Now use \scite{1.19}(1) for $\aleph_2$-c.c.c., \scite{1.16} for
properness, \scite{1.19}(2) for part (a) and \scite{1.18} for 
part (b) and \scite{1.19}(3) for part (c). \newline
${}$ \hfill$\square_{\scite{1.21}}$
\enddemo
\bigskip

\demo{\stag{1.22} Conclusion}  Let $P_{\omega_3}$ be as in \scite{1.21}.
Assuming for simplicity $V \models GCH$ in
$V^{P_{\omega_3}}$ there is no Gross space; i.e. for every countable field
$F$ and vector space ${\Cal U}$ over $F$ with inner products (i) with
dimension $\aleph_0$ for some ${\Cal U}_1,{\Cal U}_2 \subseteq {\Cal U}$,
we have dim ${\Cal U}_1 = \text{ dim } {\Cal U},\text{ dim } {\Cal U}_2 =
\aleph_0$ and ${\Cal U}_1 \bot {\Cal U}_2$.
\enddemo
\bigskip

\demo{Proof}  Let $\lambda = \text{ dim}({\Cal U})$ and $\langle x_i:i 
< \lambda \rangle$ be a basis of ${\Cal U}$.
\enddemo
\bigskip

\noindent
\underbar{Case 1}:  $\lambda = \aleph_1$.

Use clause $(b)'$ of the conclusion of \scite{1.21}.
\bigskip

\noindent
\underbar{Case 2}:  $\lambda = \aleph_2$.

Use clause $(c)$ of the conclusion of \scite{1.21}.
\bigskip

\noindent
\underbar{Case 3}:  $\lambda = \aleph_3$.

Use clause $(a)$ of the conclusion of \scite{1.21} (i.e.
let $\langle c_\alpha:\alpha \in S \rangle$ be a $\clubsuit_S$, \newline
$S \subseteq \{\delta < \aleph_2:\text{cf}(\delta) = \aleph_0\}$ stationary.

Define $f_\beta:\omega_2 \rightarrow F$ be $f_\beta(i) = (x_\beta,x_i)$.
So for each $i$ for some $\delta(i) \in S,f_\beta \restriction c_\delta$ is
constant.  So for some $\delta,B_\delta = \{\gamma < \lambda:\delta(i) =
\delta\}$ has cardinality $\lambda$.  So \newline
${\Cal U}_1 = \text{ Span } \langle x_\beta:
\beta \in B_\delta \rangle, {\Cal U}_i = \text{ Span}\{x_i - x_j:i,j \in
c_\delta\}$ are as required).
\bigskip

\noindent
\underbar{Case 4}:   $\lambda = \text{ cf}(\lambda) > \aleph_3$
or just cf$(\lambda) > \aleph_2$.

The proof is as in the case $\lambda = \aleph_3$.
\bigskip

\noindent
\underbar{Case 5}:  $\lambda > \text{ cf}(\lambda) \in \{\aleph_0,\aleph_1,
\aleph_2\}$.

Cases of $\lambda$ singular we reduce it to the problem of cf$(\lambda)$ 
as for $\mu \in (\aleph_3,\lambda)$ regular for some 
$B_\mu \in [\mu]^\mu,f_\mu:\aleph_2 \rightarrow F$ we have
$(\forall \alpha \in B_\mu)(\forall j < \aleph_2)[(x_\alpha,x_j) =
f_\mu(j)]$. \newline
${}$ \hfill$\square_{\scite{1.22}}$
\bigskip

\remark{\stag{1.23} Remark}  1) In \scite{1.21} - \scite{1.22} we can 
replace the choice $(\kappa,\mu,\alpha) = (\aleph_1,\aleph_2,\aleph_3)$ 
by any other satisfying the requirements in \scite{1.20} (e.g. $2^{< \kappa} 
= \kappa < \mu = \text{ cf}(\mu),(\forall \beta < \mu)[|\beta|^{< \kappa} 
< \mu]\text{ cf}(\alpha) > \mu,\alpha$ divisible by $|\alpha|,|\alpha| 
= |\alpha| < \mu$. \newline
2) In creature $u$ (i.e. in Definition \scite{1.5}) we can also incorporate
the information ``in $Q_\gamma$ we just want to satisfy $\mu$ conditions
among a family of $\mu$", or we consider $< \mu$ conditions to begin with.
Also making the dividing line ``$|a_\gamma| < \kappa$" or ``$a_\gamma$
countable" does not make much difference.  \newline
3) The use of $a_\gamma$ of cardinality $< \kappa$ but non-empty is when we
would like in the definition of $\theta_\gamma$ to use as a parameter a real
not in $V^{P_\gamma} \backslash V$ (or even a bounded subset of $\kappa$).
\endremark
\newpage

\shlhetal

\newpage
    
REFERENCES.  
\bibliographystyle{lit-plain}
\bibliography{lista,listb,listx}

\enddocument

\bye